\documentclass[12pt]{article}
\usepackage{amsfonts,amsmath}

\def \beq {\begin{eqnarray}}
\def \eeq {\end{eqnarray}}
\def \beqn {\begin{eqnarray*}}
\def \eeqn {\end{eqnarray*}}

\newcommand{\halmos}{\rule{1ex}{1.4ex}}

\newcounter{for}[section]

\numberwithin{equation}{section}

\newtheorem{itlemma}{Lemma}[section]
\newtheorem{itproposition}[itlemma]{Proposition}
\newtheorem{theorem}[itlemma]{Theorem}
\newtheorem{itcorollary}[itlemma]{Corollary}
\newtheorem{itremark}[itlemma]{Remark}
\newtheorem{itremarks}[itlemma]{Remarks}
\newtheorem{itdefinition}[itlemma]{Definition}
\newtheorem{itexample}[itlemma]{Example}

\newenvironment{fact}{\begin{itfact}\rm}{\end{itfact}}
\newenvironment{claim}{\begin{itclaim}\rm}{\end{itclaim}}
\newenvironment{lemma}{\begin{itlemma}}{\end{itlemma}}
\newenvironment{remark}{\begin{itremark}\rm}{\end{itremark}}
\newenvironment{remarks}{\begin{itremarks} \rm}{\end{itremarks}}
\newenvironment{corollary}{\begin{itcorollary}}{\end{itcorollary}}
\newenvironment{proposition}{\begin{itproposition}}{\end{itproposition}}
\newenvironment{definition}{\begin{itdefinition}\rm}{\end{itdefinition}}
\newenvironment{example}{\begin{itexample}\rm}{\end{itexample}}
\newenvironment{proof}{\noindent {\em Proof}.\ \
}{\hspace*{\fill}$\halmos$\medskip}
\newcommand{\be}[1]{\addtocounter{for}{1} \begin{equation}\label{#1}}
\newcommand{\ee}{\end{equation}}
\newcommand{\bl}[1]{\begin{lemma}\label{#1}}
\newcommand{\br}[1]{\begin{remark}\label{#1}}
\newcommand{\brs}[1]{\begin{remarks}\label{#1}}
\newcommand{\bt}[1]{\begin{theorem}\label{#1}}
\newcommand{\bd}[1]{\begin{definition}\label{#1}}
\newcommand{\bp}[1]{\begin{proposition}\label{#1}}
\newcommand{\bc}[1]{\begin{corollary}\label{#1}}
\newcommand{\bfact}[1]{\begin{fact}\label{#1}}
\newcommand{\bex}[1]{\begin{example}\label{#1}}
\newcommand{\ec}{\end{corollary}}
\newcommand{\efact}{\end{fact}}
\newcommand{\eex}{\end{example}}
\newcommand{\el}{\end{lemma}}
\newcommand{\er}{\end{remark}}
\newcommand{\ers}{\end{remarks}}
\newcommand{\et}{\end{theorem}}
\newcommand{\ed}{\end{definition}}
\newcommand{\ep}{\end{proposition}}
\newcommand{\epr}{\end{proof}}
\newcommand{\bpr}{\begin{proof}}
\newcommand{\bcl}[1]{\begin{claim}\label{#1}}
\newcommand{\ecl}{\end{claim}}

\newcommand{\ecs}{\end{corollary}}
\newcommand{\eers}{\end{exercise}}
\newcommand{\eexs}{\end{example}}
\newcommand{\eems}{\end{example}}
\newcommand{\els}{\end{lemma}}
\newcommand{\eles}{\end{lemmaex}}
\newcommand{\ets}{\end{theorem}}
\newcommand{\eds}{\end{definition}}
\newcommand{\eps}{\end{proposition}}

\newcommand{\bi}{\begin{itemize}}
\newcommand{\ei}{\end{itemize}}
\newcommand{\ben}{\begin{enumerate}}
\newcommand{\een}{\end{enumerate}}

\def\vbar{\mathchoice{\vrule height6.3ptdepth-.5ptwidth.8pt\kern-.8pt}
   {\vrule height6.3ptdepth-.5ptwidth.8pt\kern-.8pt}
   {\vrule height4.1ptdepth-.35ptwidth.6pt\kern-.6pt}
   {\vrule height3.1ptdepth-.25ptwidth.5pt\kern-.5pt}}
\def\fudge{\mathchoice{}{}{\mkern.5mu}{\mkern.8mu}}
\def\bbc#1#2{{\rm \mkern#2mu\vbar\mkern-#2mu#1}}
\def\bbb#1{{\rm I\mkern-3.5mu #1}}
\def\bba#1#2{{\rm #1\mkern-#2mu\fudge #1}}
\def\bb#1{{\count4=`#1 \advance\count4by-64 \ifcase\count4\or\bba A{11.5}\or
   \bbb B\or\bbc C{5}\or\bbb D\or\bbb E\or\bbb F \or\bbc G{5}\or\bbb H\or
   \bbb I\or\bbc J{3}\or\bbb K\or\bbb L \or\bbb M\or\bbb N\or\bbc O{5} \or
   \bbb P\or\bbc Q{5}\or\bbb R\or\bbc S{4.2}\or\bba T{10.5}\or\bbc U{5}\or
   \bba V{12}\or\bba W{16.5}\or\bba X{11}\or\bba Y{11.7}\or\bba Z{7.5}\fi}}

\def \T {{\cal{T}}}
\def \I {{\mathsf{I}}}

\def \A {{\cal{A}}}

\def \FF {{\cal{F}}}
\def \PP {{\cal{P}}}
\def \GG {{\cal{G}}}
\def \C {{\cal{C}}}
\def \D {{\cal{D}}}

\def \S {{\cal{S}}}

\def\sqr#1#2{{\vcenter{\vbox{\hrule height .#2pt
                             \hbox{\vrule width .#2pt height#1pt \kern#1pt
                                   \vrule width .#2pt}
                             \hrule height .#2pt}}}}

\def\pmb#1{\setbox0=\hbox{#1}%
   \kern-.025em\copy0\kern-\wd0
   \kern.05em\copy0\kern-\wd0
   \kern-.025em\raise.0433em\box0 }
\def\sqr#1#2{{\vcenter{\vbox{\hrule height.#2pt
     \hbox{\vrule width.#2pt height#1pt \kern#1pt
   \vrule width.#2pt}\hrule height.#2pt}}}}

\def\B{{\mathbb B}}
\def\N{{\mathbb N}}   
\def\Z{{\mathbb Z}}
\def\R{{\mathbb R}}

\def\P{{\mathbb P}}
\def\E{{\mathbb E}}  
\def\id{{\mathbf{1}}}

\def\p{\partial}

\def\bs{\backslash}

\def\var{{\rm var}}

\def\AA{{\mathbb A}}
\def\reff#1{(\ref{#1})}

\newcommand {\cro}[1] {\left[ {#1} \right]}
\newcommand {\acc}[1] {\left\{ {#1} \right\}}
\newcommand {\pare}[1] {\left( {#1} \right)}

\newcommand {\sous}[1] {\underline{#1}}
\newcommand{\norm}[1]{\left\| #1\right\|}

\begin{document}
\title{
   A note on fluctuations\\
   for internal diffusion limited aggregation\footnote{
    Supported by: GDRE 224 GREFI-MEFI,
    the French Ministry of Education 
    through the ANR BLAN07-2184264 grant,
    by the European Research Council 
    through the ``Advanced Grant''  PTRELSS 228032.
  }
}
\author{
  \renewcommand{\thefootnote}{\arabic{footnote}}
  Amine Asselah\footnotemark[1]
  \and
  \renewcommand{\thefootnote}{\arabic{footnote}}
  Alexandre Gaudilli\`ere\footnotemark[2]
}
\date{}

\footnotetext[1]{
    LAMA, Universit\'e Paris-Est -- e-mail: amine.asselah@univ-paris12.fr
  }
\footnotetext[2]{
  LATP, Universit\'e de Provence, CNRS,
  39 rue F. Joliot Curie,
  13013 Marseille, France\\
  \indent\indent e-mail: gaudilli@cmi.univ-mrs.fr
}
\maketitle
\markboth
    {Fluctuations for internal DLA}
    {Fluctuations for internal DLA}
\pagestyle{myheadings}
\begin{abstract}
We consider a cluster growth model on $\Z^d$, called
internal diffusion limited aggregation (internal DLA).
In this model, random walks start at the origin, one at a time,
and stop moving when reaching a site not occupied
by previous walks.
It is known that the asymptotic shape of the cluster is
spherical. Also, when dimension is 2 or more, and when the cluster
has volume $n^d$, it is known that fluctuations of the radius
are at most of order $n^{1/3}$. We improve this estimate to $n^{1/(d+1)}$,
in dimension 3 or more. In so doing, we introduce a closely
related cluster growth
model, that we call the flashing process, whose fluctuations are
controlled easily and accurately. This process is coupled to internal DLA
to yield the desired bound.
Part of our proof adapts the approach of Lawler, Bramson and Griffeath,
on another space scale, and uses a sharp estimate (written
by Blach\`ere in our Appendix) on the expected time spent
by a random walk inside an annulus.

\end{abstract}
\maketitle

\section{Introduction}
The internal DLA cluster of volume $N$, say $A(N)$,
is obtained inductively as follows. Initially, we assume
that the explored region is empty, that is $A(0)=\emptyset$.
Then, consider $N$ independent discrete-time
random walks $S_1,\dots,S_N$ starting from 0. 
Assume $A(k-1)$ is obtained, and define
\be{time-settling}
\tau_k=\inf\acc{t\ge 0:\ S_k(t)\not\in A(k-1)},\quad\text{and}\quad
A(k)=A(k-1)\cup \{S_k(\tau_k)\}.
\ee
In such a particle system, we call explorers the particles.
We say that the $k$-th explorer is {\it settled} on $S_k(\tau_k)$ after time
$\tau_k$, and is {\it unsettled} before time $\tau_k$. The cluster $A(N)$
is the positions of the $N$ settled explorers.
We study the growth of $A(N)$, as $N$ tends to infinity.

The mathematical model of
internal DLA was introduced first in the chemical physics literature
by Meakin and Deutch~\cite{meakin-deutch}.
There are many industrial processes that look like
internal DLA (see the nice review paper \cite{landolt}). 
The most important seems to be electropolishing,
defined as {\it the improvement of surface finish of a metal effected
my making it anodic in an appropriate solution}. There are actually
two distinct industrial processes (i) {\it anodic levelling or smoothing}
which corresponds to the elimination of surface roughness of
height larger than~1 micron, and (ii) {\it anodic brightening} which
refers to elimination of surface defects which are protruding by
less than~1 micron. The latter phenomenon
requires an understanding of atom removal from a crystal lattice.
It was noted in \cite{meakin-deutch}, at a qualitative level, 
that the model produces smooth clusters, and the authors wrote
``it is also of some fundamental significance to know just how smooth
a surface formed by diffusion limited processes may be''.

Diaconis and Fulton~\cite{diaconis-fulton} introduced internal DLA
in mathematics. Their model is more general than ours:
explorers can start on distinct sites, and the explored region
at time 0 is not necessarily empty.
They were interested in
defining a random growth process by iterating simple operation.
They introduced many variations,
and treat, among other things, the special one dimensional case.

In dimensions two and more, 
Lawler, Bramson and Griffeath~\cite{lawler92} prove that in order
to cover, without holes, a sphere of radius $n$, we need about the
number of sites of $\Z^d$ contained in this sphere. In other words,
the asymptotic shape of the cluster
is a sphere. Then, Lawler in~\cite{lawler95} shows
subdiffusive fluctuations.
The latter result is formulated in terms of inner and outer errors,
which we now introduce with some notation.
We denote with $\|\cdot\|$ the euclidean norm on $\R^d$.
For any $x$ in $R^d$ and $r$ in $\R$, set
\be{ball-dfn}
B(x,r) = \left\{ y\in\R^d :\: \|y-x\| < r \right\}
\quad\mbox{and}\quad \B(x,r) = B(x,r) \cap \Z^d.
\ee
For $\Lambda\subset\Z^d$,
$|\Lambda|$ denotes the number of sites in $\Lambda$.
The inner error $\delta_I(n)$ is such that
\be{def-inner}
n-\delta_I(n)=\sup\acc{r \geq 0:\: \B(0,r)\subset A(|\B(0,n)|)}.
\ee
Also, the outer error $\delta_O(n)$ is such that
\be{def-outer}
n+\delta_O(n)=\inf\acc{r \geq 0:\: A(|\B(0,n)|)\subset \B(0,r)}.
\ee
The main result of~\cite{lawler95} reads as follows.
\bt{theo-lawler}[Lawler]
Assume $d\ge 2$. With probability 1, 
\be{lawler-results}
\lim_{n\to\infty} \frac{\delta_I(n)}{n^{1/3}\log(n)^2}=0,
\quad\text{and}\quad
\lim_{n\to\infty} \frac{\delta_O(n)}{n^{1/3}\log(n)^4}=0.
\ee
\et   
Since Lawler's paper, published 15 years ago,
no improvement of these estimates was achi\-eved,
but it is believed that fluctuations are on a much smaller scale
than $n^{1/3}$. 
Computer simulations~\cite{machta-moore,levine}
suggest indeed that fluctuations are logarithmic.
In addition,
Levine and Peres studied
a deterministic analogue of internal DLA,
the rotor-router model,
introduced by J.Propp \cite{kleber-propp}. 
They bound, in \cite{levine-peres}, the inner error $\delta_I(n)$
by $\log(n)$, and the outer error $\delta_O(n)$ by $n^{1-1/d}$.

We present here an improvement on \reff{lawler-results}
in dimension $d\geq 3$.
\bt{theo-main}
Assume $d\ge 3$. There is a positive constant $A_d$ such that,
with probability 1,
\be{ag-results}
\limsup_{n\to\infty} \frac{\delta_I(n)}{n^{\frac{1}{d+1}}\log(n)}\le A_d
\quad\text{and}\quad
\lim_{n\to\infty} \frac{\delta_O(n)}{n^{\frac{1}{d+1}}\log^2(n)}= 0.
\ee
\et

Let us now recall the approach of Lawler, Bramson
and Griffeath in~\cite{lawler92}, explain our idea, and 
introduce a new growth model.
The approach of \cite{lawler92} is based on estimating
the number $W(z)$ of explorers that visit each site $z\in \Z^d$. 
It is based on the following observations.
(i) If explorers would not settle,
they would just be independent random walks; (ii) exactly one
explorer occupies each site of the cluster. Thus, if we
launch one explorer from each site of the cluster $A(N)$, and
call $M(A(N),z)$ the number of crossings of site $z$, then 
$M(N\delta_0,z)$ which we define
as $W(z)+ M(A(N),z)$ would be equal in law to the number of 
walks crossing $z$ out of $N$ independent walks started on 0.
 Even though $M(A(N),z)$ and $M(N\delta_0,z)$
are dependent variables, some estimates on $P(W(z)=0)$ can be
extracted from estimating the means of 
$M(A(N),z)$ and of $M(N\delta_0,z)$.

Rather than thinking in terms of one single
site, we observe that a site has good chances 
to lie inside the cluster if
some {\it large} region, about this site,
is crossed by {\it many} explorers.
How {\it large} should be this region,
say $\C$, and how {\it many} should be these
crossings, say $W(\C)$, can be answered as follows.
The size and location of $\C$ should be such that
(i) the expected number of crossings of $\C$ is much larger than
its standard deviation, and (ii) the number of crossings 
of $\C$ needed to cover $\C$
is of order $|\C|$ 
(as suggested by the spherical shape result of~\cite{lawler92})
\be{main-heuristics}
\text{(i)}\quad E[W(\C)]\gg\sqrt{\var(W(\C))},\qquad\text{(ii)}\quad
E[W(\C)] \geq |\C|.
\ee
Now, assume that $|\B(0,n)|$ explorers start at the origin.
For a space-scale $h(n)$ and an integer $k>1$, to be determined,
assume that $\B(0,n-kh(n))$ is covered by settled explorers.

Partition the shell $\S=\B(0,n-(k-1)h(n))\bs \B(0,n-kh(n))$
into about $(n/h(n))^{d-1}$ {\it cells}, each of volume $h(n)^d$.
Cells are brick-like domain, of side length the width of the shell.
It is convenient to imagine that a cell $\C\in \S$ is the basis
of a column of $k$ cells reaching the boundary of $\B(0,n)$. It is
also convenient to stop the explorers as they reach the boundary
of $\B(0,n-kh(n))$. Thus, with such a stopped process, explorers are
either settled inside the $\B(0,n-kh(n))$ or
unsettled but stopped on its boundary,
that we denote by $\p\B(0,n-kh(n))$. What we have called earlier
{\it the number of explorers crossing $\C$} is taken here
to be the unsettled explorers
stopped on $\C\cap \p\B(0,n-kh(n))$. In these heuristics, we make
the simplifying assumption that only explorers stopped on
$\C\cap \p\B(0,n-kh(n))$ can cover $\C$ once released.
\begin{itemize}
\item
On the average, $k h(n)^d$ explorers are stopped on each cell of $\S$.
Thus, for $k$ large enough
(ii) of \reff{main-heuristics} should be fulfilled.
\item
The standard deviation of $W(\C)$ is a delicate issue not yet settled.
Here follows a heuristic justification
for an upper bound for the standard deviation.
If each explorer had to choose uniformly at random
a cell ${\cal C}$ of $\S$
and were constrained to perform a reflected random walk inside a cone
issued from the origin and with base $\C$, then we would
obtain another growth model whose fluctuations are expectedly
larger than those of internal DLA.
For such a model the standard deviation of the crossing of a cell
is of order   
\be{heur-alex}
\sqrt{|\B(0,n)| P_0(S(\tau)\in\C)}\sim
\sqrt{n^d\left(h(n)/n\right)^{d-1}},
\ee
where $\tau$ is the exit time from $\B(0,n-kh(n))$, and $P_0$ is the law of
a simple random walk, $S$, starting at the origin.
We expect \reff{heur-alex}
to be an upper bound for the standard deviation of $W(\C)$.
\end{itemize}
At a heuristic level, (i) of \reff{main-heuristics} follows if
\be{heuristics-1}
k h(n)^d \gg\sqrt{n^d \times
\left(\frac{h(n)}{n}\right)^{d-1}}
\Longleftarrow h(n) \gg n^{\frac{1}{1+d}}.
\ee

This discussion motivates 
a growth model associated from the start
with the exponent $\frac{1}{1+d}$.
We call this model the flashing internal DLA process,
or simply {\it the flashing process}.
The flashing process looks like internal DLA on a large scale, but
has a distinct covering mechanism, which makes it much simpler to analyze.
To obtain this growth model, we generalize the rules described in 
\reff{time-settling}, by enabling explorers to
settle only at special times, called flashing times.
Thus, each explorer is associated
with a sequence of stopping times, and it is only at these times
that it settles if outside the cluster.
The precise definition of the chosen stopping times requires
additional notation, which we postpone
to Section~\ref{sec-flashing}. We describe here
key features of the flashing process.

First, $\Z^d$ is partitioned into concentric
shells around the origin: a shell at a distance $r$ from
the origin has a width of order $r^{\frac{1}{d+1}}$.
Each shell is in turn partitioned into cells, which are brick-like domain,
of side length the width of the shell.
The key features are as follows.
\begin{itemize}
\item[$\PP$a.] An explorer {\it flashes} at most once in each shell.
\item[$\PP$b.]  The {\it flashing} position in a shell,
is essentially uniform
over the cell an explorer first hits upon entering the shell.
\item[$\PP$c.]  When an explorer leaves a shell, it cannot afterward
{\it flash} in it.
\end{itemize}
Feature $\PP$b is the seed of a deep difference with internal DLA.
{\it The mechanism of covering a cell, for the flashing process, 
is very much the same
as completing an album in the classical coupon-collector process}.
Thus, we need of the order of $V\log(V)$ explorers 
to cover a cell of volume $V$. For internal DLA, with explorers
started at the origin, we need only of order $V$ explorers to
cover a sphere of volume $V$ as shown in \cite{lawler92}, and 
we believe that we need a number of explorers
of order $|{\cal C}|$ to cover a cell ${\cal C}$.

Feature $\PP$c is essential for having the following coupling
between flashing and internal DLA processes.
\bt{theo-coupling}
There is a coupling between the two processes such that,
for all $k,N\geq 1$, and for a sequence $\{r_k,k\in \N\}$
going to infinity
(that we describe in Section~\ref{sec-flashing}), 
and $h_k=r_k^{1/(d+1)}$ 
\begin{itemize}
\item 
  if $A^*(N)\subset \B(0,r_k + h_k )$,
  then $A(N)\subset \B(0,r_k + h_k )$,
\item 
  if $\B(0,r_k + h_k )\subset A^*(N)$, 
  then $\B(0,r_k + h_k )\subset A(N)$.
\end{itemize}
\et
For the flashing process, 
we control easily the inner error.
Then, to control the outer error we follow the approach of~\cite{lawler95},
though with a simpler proof
(allowed by the cell structure used to build our growth model).

\bt{theo-flashing}
There is a positive constant $A_d$ such that, \reff{ag-results}
hold for the flashing process.
\et
Also, we know that flashing internal DLA does exhibit
power-law fluctuations.
\bt{theo-optimal}
With probability 1,
\be{main-prop2}
\lim_{n\to\infty} \frac{n^{\frac{1}{d+1}}}{\delta^*_I(n)}=0.
\ee
\et
Theorem~\ref{theo-coupling} and Theorem~\ref{theo-flashing} imply
Theorem~\ref{theo-main}.

Let us now describe the heuristics behind Theorem~\ref{theo-optimal}.
It is useful to organize the flow of explorers in the flashing
process into {\it exploration waves}, in the way of Section 3
of \cite{lawler95}. That is, in the $k$-th exploration wave,
the explorers either stop as they reach
the {\it bulk} of $\S_k$, or settle before reaching $\S_k$.
Consider now the exploration wave associated with the last shell
making $\B(0,n)$, say $\S^*$.
Assume that at this time, the cluster fills $\B(0,n)\bs \S^*$.
From our definition, $\S^*$ has width of order $n^{\frac{1}{d+1}}$, and
the last shell receives a number of stopped explorers
equal to its volume. There is necessarily one cell in $\S^*$ which
receives of the order of its own volume, and for a coupon-collector
process, it is very unlikely that the explorers stopped  in this very cell
can cover the bulk of this cell before escaping $\S^*$. By feature
$\PP$c, if a hole is left in $\S^*$ after the explorers leave $\S^*$,
this hole remains uncovered forever.
These heuristics pose two questions concerning internal DLA,
we are unable to answer at the moment.
\begin{itemize}
\item How many explorers stopped in {\it the bulk} of
a cell are needed to cover the whole cell?
\item What is the correct order of fluctuations in internal DLA?
\end{itemize}

The rest of the paper is organized as follows.
Section~\ref{sec-notation} introduces the main notation, and
recall well known useful facts.
In Section~\ref{sec-flashing}, we build the flashing process,
give an alternative construction, 
and prove Theorem~\ref{theo-coupling}. In Section~\ref{sec-harmonic},
we obtain a sharp estimate on the expected number of explorers
crossing a given cell, and prove $\PP$-b. Both proofs
are based on classical potential theory estimates.
In Section~\ref{sec-results}, we prove Theorems~\ref{theo-flashing}
and \ref{theo-optimal}. Finally, in the Appendix,
S\'ebastien Blach\`ere gives a sharp estimate
on the expected time spent in an annulus by a random walk.
\section{Notation and useful tools}\label{sec-notation}
\subsection{Notation}
We say that $z,z'\in \Z^d$ are nearest neighbors when $\|z-z'\|=1$, and
we write $z\sim z'$. For any subset $\Lambda\subset \Z^d$, we define
\be{not2}
\p \Lambda=\acc{z\in \Z^d\bs \Lambda:\: \exists z'\in \Lambda, z'\sim z}.
\ee
For any $r\leq R$ we define the annulus
\be{annulus-dfn}
A(r,R) = B(0,R)\setminus B(0,r)
\quad
\mbox{and}
\quad
{\mathbb A}(r,R) = {\mathbb A}(r,R)\cap\Z^d
\ee
A trajectory $\gamma$ is a discrete nearest-neighbor path on $\Z^d$.
That is $\gamma: \N\to\Z^d$ with $\gamma(t)\sim\gamma(t+1)$ for all $t$.
The law of the simple random walk started in $z$, is denoted with $\P_z$.
For a subset $\Lambda$ in $\Z^d$, and a trajectory $\gamma$,
we define the hitting time of $\Lambda$ as
\[
H(\Lambda;\gamma)=\min\{t\geq 0 :\: \gamma(t)\in\Lambda\}.
\]
We often omit $\gamma$ in the notation when no confusion is possible.
We use the shorthand notation
\[
B_n=B(0,n),\quad \B_n=\B(0,n), \quad H_R=H(B_R^c),
\quad \text{and}\quad H_z = H(\{z\}).
\]
For any $a$,$b$ in $\R$ we write $a\wedge b = \min\{a,b\}$, and
$a\vee b=\max\{a,b\}$.
Let $\Gamma$ be a finite collection of trajectories
on $\Z^d$. For $R>0$, $z$ in $\Z^d$ and $\Lambda$ a subset
of $\Z^d$, we call
$M(\Gamma, R,z)$ (resp. $M(\Gamma,R, \Lambda)$)
the number of trajectories which exit $\B(0,R)$ on $z$ (resp. in $\Lambda$):
\be{def-M}
M(\Gamma,R,z) = \sum_{\gamma\in \Gamma} \id_{\{\gamma(H_R)=z\}},
\quad\text{and}\quad
M(\Gamma,R,\Lambda) = \sum_{z\in\Lambda} M(\Gamma,R,z).
\ee
When we deal with a collection of independent random trajectories, we
rather specify its initial configuration $\eta\in \N^{\Z^d}$,
so that $M(\eta,R,z)$ is the number of random walks
starting from $\eta$ and hitting $\B(0,R)^c$ on $z$. Two types of
initial configurations are important here:
(i) the configuration $n\id_{z^*}$ formed by $n$ walkers starting on
a given site $z^*$,
(ii) for $\Lambda\subset\Z^d$,
the configuration $\id_\Lambda$
that we simply identify with $\Lambda$.
For any configuration $\eta\in\N^{\Z^d}$ we write
\be{not5}
|\eta| = \sum_{z\in\Z^d} \eta(z).
\ee

We are in dimension 3 or more, and Green's function
of the simple random walk is well defined and denoted $G$. That is,
for any $x,y\in \Z^d$
\be{def-green}
G(x,y)=\E_x\left[ \sum_{n\ge 0} \id_{\{S(n) = y\}} \right].
\ee
For any $\Lambda\subset\Z^d$, we define Green's function
restricted to $\Lambda$, $G_\Lambda$, as follows. For $x,y\in \Lambda$
\be{green-lambda}
G_\Lambda(x,y)= \E_x\left[
  \sum_{0\le n < H(\Lambda^c)} \id_{\{S(n) = y\}} \right].
\ee
\subsection{Some useful tools}\label{sec-green}
We recall here some well known facts.
Some of them are proved for the reader's convenience. This section
can be skipped at a first reading.

In~\cite{lawler92}, the authors emphasized the fact
that the spherical limiting shape of internal DLA
was intimately linked to strong isotropy properties
of Green's function. This isotropy is expressed
by the following asymptotics (Theorem 4.3.1 of \cite{lawler-limic}).
In $d\ge 3$, there is a constant $K_g$, such that for any $z\not=0$,
\be{main-green}
\big|G(0,z) - \frac{C_d}{\|z\|^{d-2}} \big|\le \frac{K_g}{\|z\|^d}
\quad
\text{with}\quad C_d=\frac{2}{v_d(d-2)},
\ee
where $v_d$ stands for the volume of
the euclidean unit ball in $\R^d$.
The first order expansion \reff{main-green}
is proved in \cite{lawler-limic}
for general symmetric walks
with finite $d+3$ moments and vanishing third
moment. All the estimates we use are eventually based on \reff{main-green}
and we emphasize the fact that the estimate is uniform in $\|z\|$.

The following lemma is also used in the Appendix.
\bl{lem-border}
Each $z^*$ in $\Z^d\setminus\{0\}$
has a nearest-neighbor $z$ (i.e. $z^*\sim z$) such that
\be{goodnn}
\|z\|\leq \|z^*\| -\frac{1}{2\sqrt{d}}.
\ee
\el
\bpr
Without loss of generality
we can assume that all the coordinates
of $z^*$ are non-negative.
Let us denote by $b$ the maximum of these coordinates and note that
\be{dfna}
\|z^*\|^2 \leq db^2,\quad\text{and}\quad b\ge 1.
\ee
Denote by $z$ the nearest-neighbor obtained from $z^*$
by decreasing by one unit a maximum coordinate.
Using \reff{dfna}
\be{simple}
\|z^*\|^2 - \|z\|^2=b^2-(b-1)^2=2b-1\ge b\ge \frac{\|z^*\|}{\sqrt{d}}.
\ee
Note that \reff{goodnn} follows from $2\|z^*\|(\|z^*\| - \|z\|)\ge
\|z^*\|^2 - \|z\|^2$, and \reff{simple}.
\epr

We recall a rough but useful result about the exit site distribution
from a sphere. This is Lemma 1.7.4 of \cite{lawler}.
\bl{lem-lawler}
There are two positive constants $c_1,c_2$ such that for
any $z\in \p B(0,n)$, and $n>0$
\be{hitting-position}
\frac{c_1}{n^{d-1}}\le \P_0(S(H_n)=z)\le \frac{c_2}{n^{d-1}}.
\ee
\el

Finally, we recall a well known large deviations
estimate for independent Bernoulli variables
(see for instance Lemma 4.3 of \cite{bramson-lebowitz}).
\bl{lem-BL} 
For any positive integer $n$, and $\{X_1,\dots,X_n\}$
$n$ independent
Bernoulli variables, 
we have for any $x>0$,
and with $X=X_1+\dots+X_n$
\be{eq-BL}
\max\pare{P\pare{ X-E[X]\ge x}, P\pare{ X-E[X]\le -x}}
\le \exp\pare{-\pare{\min\pare{\frac{x^2}{4\var[X]},\frac{x}{2}}}}.
\ee
\el
\br{rem-BL}
Note that for a sum of Bernoulli, $\var[X]\le E[X]$, and the following
inequality is useful
\be{eq-BL1}
\max\pare{P\pare{ X-E[X]\ge x}, P\pare{ X-E[X]\le -x}}
\le \exp\pare{-\pare{\min\pare{\frac{x^2}{4E[X]},\frac{x}{2}}}}.
\ee
Also, note that if $E[X]\le E[Y]$, where $Y$ is a sum of
$m$ independent Bernoulli variables, then
\be{eq-BL2}
\max\pare{P\pare{ X-E[X]\ge x}, P\pare{ Y-E[Y]\le -x}}
\le \exp\pare{-\pare{\min\pare{\frac{x^2}{4E[Y]},\frac{x}{2}}}}.
\ee
\er
%

\section{The flashing process}\label{sec-flashing}
In this section, we construct the flashing process.
We then present a useful alternative construction of the same process.
Finally, we prove Theorem~\ref{theo-coupling} which couples the two
processes.

\subsection{Construction of the process}

\paragraph{Partitioning the lattice.}
We partition the lattice
into shells $({\cal S}_j :\: j \geq 0)$.
For a given parameter $h_0 > 0$
the first shell ${\cal S}_0$ is the ball ${\mathbb B}(0,h_0)$.
The next shells are the annuli
\begin{equation}
{\cal S}_j = {\mathbb A}(r_j - h_j, r_j + h_j),
\qquad j\geq1,
\end{equation}
where $r_j$ and $h_j$ are defined inductively by
$r_1 - h_1 = h_0$, and for $j\ge 1$
\be{dfn-r-h}
r_{j+1} - h_{j+1} = r_j + h_j, \quad\text{and}\quad 
h_j = r_j^{\frac{1}{d+1}}. 
\ee
We omit the easy check that \reff{dfn-r-h} yields
\begin{equation}
r_j \sim \left(\frac{2d}{d+1} j\right)^{\frac{d+1}{d}}.
\end{equation}
We also define
\begin{equation}
\Sigma_0 =\{0\}
\quad\mbox{and}\quad
\Sigma_j = \partial {\mathbb B}(0,r_j),
\; j\geq 1.
\end{equation}

\paragraph{Flashing times.}
Consider $\{X_j,Y_j,j\ge 0\}$ a sequence of independent Bernoulli 
variables such that
\begin{equation}
P(X_j = 1) = 1 - P(X_j = 0) = \frac{1}{h_j^d},
\end{equation}
\begin{equation}
P(Y_j = 1) = 1 - P(Y_j = 0) = \left\{
  \begin{array}{ll}
    1 &\mbox{if } j = 0,\\
    \frac{1}{2} &\mbox{if } j\geq 1,
  \end{array}
\right.
\end{equation}
Consider also a sequence of continuous independent variables
$\{R_j,j\ge 0\}$ each of which has density $g_j:[0,h_j]\to \R^+$ with
\be{density}
g_j(h)=\frac{dh^{d-1}}{h_j^d}.
\ee
For $j\ge 0$, and $z_j$ in $\Sigma_j$,
let $S$ be a random walk starting in $z_j$, an define
a stopping time $\sigma$ as follows.
If $R_j = h$ for some $h \leq h_j$ then
\begin{equation}
\sigma = \left\{
  \begin{array}{ll}
    0 
    &\mbox{if $ X_j = 1$,}\\
    H({\mathbb B}(z, h\wedge (r_j + h_j - \|z_j\|))^c)
   &\mbox{if $X_j = 0$ and $Y_j = 1$,}\\
    H({\mathbb A}(r_j - h, r_j + h)^c)
    &\mbox{if $X_j = 0$ and $Y_j = 0$.}
  \end{array}
\right.
\end{equation}
We set $H_j = H(\Sigma_j)$, and we define  the stopping times
$(\sigma_j :\: j\geq 0)$ as
\begin{equation}
\sigma_j = H_j + \sigma(S\circ\theta_{H_j}),
\end{equation}
where $\theta$ stands for the usual time-shift operator.
For $j \geq 0$ we note that, by construction, $S(t)\in{\cal S}_j$
for all $t$ such that $H_j \leq t < \sigma$ and we say that
$\sigma_j$ is a {\em flashing time} when $S(\sigma_j)$
is contained in the intersection between ${\cal S}_j$
and the cone with base $B(S(H_j), h_j / 2)$.
We call such an intersection {\it a cell centered at} $S(H_j)$,
that we denote $\C(S(H_j))$.
In other words, for any $z\in \Sigma_j$
\be{def-cell}
{\cal C}(z) = {\cal S}_j \cap\left\{
  x \in {\mathbb R}^d :\:  
  \exists \lambda \geq 0, \exists y\in B(z, h_j/2), x = \lambda y
\right\}.
\ee

\paragraph{The uniform hitting property}
The main property of the hitting time $\sigma$ constructed above
is the following proposition, which yields property $\PP$b of
the flashing process to be defined soon.
\bp{prop-uniform}
There are two positive constants $\alpha_1 < \alpha_2$,
such that, for $j\geq 0$, $z_j\in\Sigma_j$, and
$z^*\in \C(z_j)$.
\be{uniform-inequality}
\frac{\alpha_1}{h_j^d} \leq \P_{z_j}\left( S(\sigma)=z^*\right)
\leq \frac{\alpha_2}{h_j^d}.
\ee
\ep
The proof of Proposition~\ref{prop-uniform} is given in
Section~\ref{sec-harmonic}.

\paragraph{The flashing process.}
Consider a family of $N$ independent random walks
$(S_i^* : 1\leq i \leq N)$
with their hitting times, 
and stopping times $(H_{i,j}, z_{i,j}, \sigma_{i,j} :\: j\geq 0)$.
Let also $z_{i,j}=S_i(H_{i,j})$ be the first hitting position
of $E_i$ on $\Sigma_i$.

We define the cluster inductively. Set
${\cal A}^*(0) = \emptyset$. For $i\geq 1$, 
we define $\tau_i^*$ as the first flashing time
associated with $S_i^*$ when the
explorer stands outside ${\cal A}^*(i -1)$. In other words,
\begin{equation}
\tau_i^* =\min\left\{
  \sigma_{i,j}:\ j \geq 0,\
  S_i^*(\sigma_{i,j}) \in {\cal C}(z_{i,j})\cap{\cal A}^*(i-1)^c
\right\},
\end{equation}
and
\begin{equation}
  {\cal A}^*(i) = {\cal A}^*(i-1)\cup \left\{S_i^*(\tau_i^*)\right\}.
\end{equation}

\subsection{Exploration Waves}\label{sec-wave}

Rather than building ${\cal A}^*(N)$
following the whole journey of one explorer
after another, we can build ${\cal A}^*(N)$
as an increasing union of clusters formed
by stopping explorers on successive shells. 
Similar wave constructions are introduced
in~\cite{lawler92} and~\cite{lawler95},
with an equality in law between alternative constructions.
However, the features of the flashing
process are such that in our case 
the two constructions are strictly equivalent.
We use this alternative construction
in the proof of Theorem~\ref{theo-coupling}. 

We denote by $\xi_k\in({\mathbb Z}^d)^N$
the explorers positions after the $k$-th wave.
We denote by ${\cal A}^*_k(N)$
and the set of sites where settled explorers are
after the $k$-th wave. Our construction will be such that
\begin{equation}
\xi_k(i)\not\in\Sigma_k
\;\Leftrightarrow\;
\xi_k(i)\in\cup_{j<k}{\cal S}_j
\;\Leftrightarrow\;
\xi_k(i)\in {\cal A}^*_k(N).
\end{equation}
For $k = 0$ we set $\xi_0(i)= 0$, and
${\cal A}^*_0(i) = \emptyset$, for $1\leq i\leq N$.
Then, for all $k\geq 0$,
we set ${\cal A}_{k+1}^*(0) = {\cal A}^*_{k}(N)$.
For $i$ in $\{1,\cdots,N\}$, we set the following.
\begin{itemize}
\item If $\xi_k(i)\not\in\Sigma_k$, then
\[
\xi_{k+1}(i) = \xi_{k}(i) \in \cup_{j<k}{\cal S}_j,\quad\text{and}\quad
{\cal A}^*_{k+1}(i) = {\cal A}^*_{k+1}(i-1).
\]
\item
If $\xi_k(i)\in\Sigma_k$
  and $S_i(\sigma_{i,k})\in{\cal C}(z_{i,k})\cap{\cal A}^*_k(i-1)^c$, then
\[
\xi_{k+1}(i) = S_i(\sigma_{i,k})\in{\cal S}_k,\quad\text{and}\quad
{\cal A}^*_{k+1}(i) = {\cal A}^*_{k+1}(i-1) 
\cup \left\{S_i(\sigma_{i,k})\right\}.
\]
\item
If $\xi_k(i)\in\Sigma_k$
  and $S_i(\sigma_{i,k})\not\in{\cal C}(z_{i,k})\cap{\cal A}^*_k(i-1)^c$, then
\[
\xi_{k+1}(i) = S_i(H_{i,k+1})\in\Sigma_{k+1},\quad\text{and}\quad
{\cal A}^*_{k+1}(i) = {\cal A}^*_{k+1}(i-1).
\]
\end{itemize} 
In words, for each $k\geq 1$, during the $k$-th wave of exploration,
the unsettled explorers move one after the other
in the order of their labels until either settling in ${\cal S}_{k-1}$,
or reaching $\Sigma_{k}$ where they stop.
We then define ${\cal A}^*(N)$ by
\begin{equation}
{\cal A^*}(N) = \bigcup_{k\geq 1}{\cal A}^*_k(N).
\end{equation}
We explain now why this construction yields the same
cluster as our previous definition.
An explorer cannot settle inside a shell it has left,
and thus cannot settle in any shell ${\cal S}_j$
with $j<k$ if it reaches ${\Sigma}_k$. Now, since
each wave of exploration is organized according to the label ordering,
the fact that an explorer has to wait for the following explorers
before proceeding its journey beyond $\Sigma_k$ does not interfere
with the site where it eventually settles.
 
\subsection{Coupling internal DLA and flashing processes}

We use here the first definition of the flashing process, and
realize the internal DLA process using the same randomness.
\bp{prop-coupling} 
There is a coupling between
the flashing and original internal DLA processes such that,
for all $N\geq 1$,
\begin{equation}
\label{orbites}
\bigcup_{i=1}^N 
\left\{
  S_i(t) :\: 0 \leq t \leq \tau_i
\right\}
\;\subset\;
\bigcup_{i=1}^N 
\left\{
S_i^*(t) :\: 0 \leq t \leq \tau^*_i
\right\}
\end{equation}
and there is a one to one mapping
$\psi_N : {\cal A}(N) \rightarrow {\cal A}^*(N)$ such that
for $z\in \A(N)$
\begin{equation}
\text{if for }\quad k\geq 1,\quad 
z\not\in\bigcup_{j<k}{\cal S}_j,\quad\text{then}\quad
\psi_N(z)\not\in\bigcup_{j<k}{\cal S}_j.
\label{stable}
\end{equation}
\ep
Theorem~\ref{theo-coupling}
is a simple consequence of Proposition~\ref{prop-coupling}.
On the one hand,
if ${\cal A}^*(N)\subset\cup_{j<k}{\cal S}_j$
for some $k\geq 1$,
then, recalling that a flashing explorer
cannot settle a shell it has left,
the orbits of the $N$ flashing
explorers are all contained in $\cup_{j<k}{\cal S}_j$
and, by~(\ref{orbites}), so is ${\cal A}(N)$.
On the other hand,
if $\cup_{j<k}{\cal S}_j\subset{\cal A}^*(N)$,
then,~(\ref{stable}) implying that with such a coupling
\begin{equation}
\left|{\cal A}(N)\cap\cup_{j<k}{\cal S}_j\right| 
\geq 
\left|{\cal A}^*(N)\cap\cup_{j<k}{\cal S}_j\right| 
= \left|\cup_{j<k}{\cal S}_j\right|, 
\end{equation}
which implies that $\cup_{j<k}{\cal S}_j\subset{\cal A}(N)$.

\paragraph{Proof of Proposition~\ref{prop-coupling}.}
We build the coupling together with the map $\psi_N$ by induction on $N$.
We use the trajectories of the flashing
explorers to drive the internal DLA trajectories. 
We need a little more notation to do so.
For each $i\leq N$, set for simplicity $g^*(i)=S_i^*(\tau_i^*)$,
and denote by $t_{i,N}$ the length of the flashing trajectory
$(S_i^*(t) :\: 0\leq t\leq t_{i,N})$
used to form the trajectories of the original explorers.
Necessarily, we need $t_{i,N}\leq \tau_i^*$
to have~(\ref{orbites}).
For convenience, we partition ${\cal A}(N)$ into 
blue sites, say ${\cal B}(N)$, and red sites, say ${\cal R}(N)$.

We build a one to one map
$f_N : {\cal A}(N) \rightarrow \{1,\cdots,N\}$, 
together with the blue-red partition 
as follows. For each $z$ in ${\cal A}(N)$, there are two possibilities.
Either there is $i\leq N$ such that $S_i^*(\tau_i^*)=z$
and $t_{i,N} = \tau_i^*$,
and we say that $z\in {\cal B}(N)$ and we set $f_N(z)=i$.
Otherwise $z\in {\cal R}(N)$. We then define $f_N(z)$ as the label~$i$
  of the flashing explorer that was driving the random
  walk $S_j$ when the $j$-th explorer settled 
  in $z$, and this will imply, by induction, that $t_{i,N} < \tau_i$.

Finally the one to one map $\psi_N: {\cal A}(N) \to {\cal A}^*(N)$ 
is the composition $g^*\circ f_N$. Note, first, that
for all $z\in {\cal B}(N)$, $\psi_N(z) = z$,
second, 
${\cal B}(N) \subset {\cal A}(N) \cap {\cal A}^*(N)$
and last,
\[
z \in {\cal B}(N) \;\Leftrightarrow\; t_{f_N(z), N} = \tau^*_{f_N(z)} ,
\qquad
z \in {\cal R}(N) \;\Leftrightarrow\; t_{f_N(z), N} < \tau^*_{f_N(z)}.
\]
We now start our induction with $N=1$.
The first trajectory $(S_1^*(t) :\: 0 \leq t \leq \tau_1^*)$
ends in $g^*(1)$.
We use this trajectory to build that of the first internal DLA-explorer.
We set $S_1(0) = S_1^*(0)$ and immediately stop here since the origin
$0 = S_1(0)$ was initially unoccupied. As a consequence 
$t_{1,1} = 0$, $\tau_1 = 0$
and there are two possibilities:
either $\tau_1^* = 0$, ${\cal B}(1) =\{0\}$, ${\cal R}(1) = \emptyset$
and $f_1(0) = 1$, or, $\tau_1^* > 0$,
${\cal R}(1) =\{0\}$, ${\cal B}(1) = \emptyset$ and $f_1(0) = 1$.

Assume now that we have built from the trajectories  
$\{(S_i^*(t) :\: 0 \leq t \leq \tau_i^*),\ i\leq N\}$,
the clusters ${\cal A}^*(N)$, and the sets
${\cal B}(N)$ and ${\cal R}(N)$,
together with the times $\{t_{i,N},\ i\leq N\}$,
and the one to one map $f_N : {\cal A}(N) \rightarrow \{1,\cdots,N\}$.
We launch a new flashing explorer with
trajectory
$(S_{N+1}^*(t) :\: 0 \leq t \leq \tau_{N+1}^*)$
that ends in $g^*(N+1)$, and we start to define
the $(N+1)$-th trajectory for the original internal DLA process
by following $S_{N+1}^*$:
\begin{equation}
  S_{N+1}(0) = S^*_{N+1}(0),\;
  S_{N+1}(1) = S^*_{N+1}(1),\;
  \dots
\end{equation}
\begin{itemize}
\item
  If $\{S_{N+1}^*(t) :\: 0 \leq t \leq \tau_{N+1}^*\}$
  is not contained in ${\cal A}(N)$
  then $S_{N+1}$ settles the first time $S_{N+1}^*$ exits ${\cal A}(N)$,
  that is at time (resp. on a site $z$)
\[
t_{N+1,N+1}=\inf\{k\ge 0: S^*_{N+1}\not\in A(N)\},\quad
\pare{\text{resp.}\quad
z = S^*_{N+1}(t_{N+1,N+1})\in {\cal A}(N)^c}.
\]
  We then set 
\[
\forall i\leq N,\quad t_{i,N+1} = t_{i,N}, \quad
f_{N+1}|_{{\cal A}(N)} = f_N,\quad\text{and}\quad 
f_{N+1} (z) = N+1,
\]
and, if $t_{N+1,N+1}= \tau_{N+1}^*$ (resp. $t_{N+1,N+1}< \tau_{N+1}^*$).
\[
\begin{split}
{\cal B}(N+1) =& \{z\}\cup {\cal B}(N)\quad
(\mbox{resp. } {\cal B}(N)),\\
{\cal R}(N+1) =& {\cal R}(N)\quad
(\mbox{resp. } \{z\} \cup {\cal R}(N)),
\end{split}
\]
\item
  If $\{S_{N+1}^*(t) :\: 0 \leq t \leq \tau_{N+1}^*\}$
  is contained in ${\cal A}(N) = {\cal B}(N) \cup {\cal R}(N)$
  then $S_{N+1}^*$ settles necessarily 
  in a red site $z$ since ${\cal B}(N) \subset {\cal A}^*(N)$.
  This red site is occupied by an explorer
  that was driven by a flashing explorer $i = f_N(z)$
  when it settled, and we have
  $t_{i,N} < \tau_i^*$.
  After reaching $S_{N+1}(\tau^*_{N+1}) = S^*_{N+1}(\tau^*_{N+1})$
  in our definition of $S_{N+1}$,
  we set
  \begin{equation}
    S_{N+1}(1+\tau^*_{N+1}) = S^*_{i}(1+t_{i,N}),\;
    S_{N+1}(2+\tau^*_{N+1}) = S^*_{i}(2+t_{i,N}),
    \dots
  \end{equation}
  we set $t_{N+1, N+1} = \tau^*_{N+1}$,
  we turn blue the site $z$, we set $f_{N+1}(z) = N+1$
  and we proceed as previously:
    If $\{S_i^*(t) :\: t_{i,N} < t \leq \tau_i^*\}\not\subset {\cal A}(N)$
    then $S_{N+1}$ settles at the first exit from ${\cal A}(N)$
    in some $z' = S^*_i(t_{i,N+1})\in {\cal A}(N)^c$
    with $t_{i,N+1}\leq \tau_i^*$
    and we set 
\[
\forall j\leq N,\quad j\neq i,\ t_{j,N+1} = t_{j,N},\quad 
f_{N+1}|_{{\cal A}(N)\setminus \{z\}} = f_N|_{{\cal A}(N)\setminus \{z\}},
\quad\text{and}\quad f_{N+1} (z') = i,
\]
and if $t_{i,N+1}= \tau_i^*$ (resp. $t_{i,N+1}< \tau_i^*$)
\[
\begin{split}
      {\cal B}(N+1) =& \{z'\}\cup {\cal B}(N) \cup \{z\}\quad
      (\mbox{resp. } {\cal B}(N) \cup \{z\}),\\
      {\cal R}(N+1) =& {\cal R}(N) \setminus \{z\}\quad 
      (\mbox{resp. } \{z'\} \cup {\cal R}(N) \setminus \{z\}).
\end{split}
\]
Otherwise, $\{S_i^*(t) :\: t_{i,N} < t \leq \tau_i^*\} \subset {\cal A}(N)$
    and $S_i^*$ settles on a red site $z'$. With $i' = f_N(z')$
    we have $t_{i',N} < \tau_{i'}^*$
    and after
    $S_{N+1}(\tau_i^*-t_{i,N}+\tau^*_{N+1}) = S^*_i(\tau^*_i)$
    we can set
    \begin{equation}
      S_{N+1}(1+\tau_i^*-t_{i,N}+\tau^*_{N+1})
      =S^*_{i'}(1+t_{i',N}),\;
     \dots
    \end{equation}
    we set $t_{i,N+1} = \tau^*_i$,
    we turn blue the site $z'$, we set $f_{N+1}(z') = i$ and so on.
\end{itemize}
Since the number of red sites is finite
this procedure necessarily reaches an end
and one immediately checks
that we define in this way
${\cal A}_{N+1} = {\cal B}_{N+1} \cup {\cal R}_{N+1}$
together with the times $\{t_{i,N+1},\ i\leq N+1\}$,
and a function $f_{N+1} : {\cal A}(N+1) \rightarrow \{1,\dots, N\}$
with all the required properties
($f_{N+1}$ is one to one and 
$t_{f_{N+1}(z), N+1} < \tau^*_{f_{N+1}(z)}$
for all $z$ in ${\cal R}_{N+1}$).

Note first that we have the orbits inclusion by construction.
Now, the law of $({\cal A}(N) :\: N\geq 1)$
is that of the internal DLA process.
Indeed, the part of the flashing trajectories 
$\{(S_i^*(t) :\: t_{i,N} < t \leq \tau^*_i),\ i\leq N\}$, that can be used 
together with $(S_{N+1}^*(t) :\: 0\leq t \leq \tau^*_{N+1})$,
to build ${\cal A}(N+1)$ have increments
that are independent from ${\cal A}(N)$.

Finally, for any $k\geq 1$, 
$\psi_N : {\cal A}(N) \rightarrow {\cal A}^*(N)$
does associate with any site
outside $\cup_{j<k} {\cal S}_j$
a site outside $\cup_{j < k} {\cal S}_j$.
Indeed, with any blue site, $\psi_N = g^*\circ f_N$ associates 
that site itself, while with each red site $\psi_N$ associates 
the end point of a flashing trajectory that visits that site.
And a flashing trajectory that exits $\cup_{j<k} {\cal S}_j$
necessarily settles outside $\cup_{j<k} {\cal S}_j$.

\section{Estimates on the Harmonic measure}\label{sec-harmonic}
We gather in this section two results which deal with
the hitting probability of sets. The first one relies
on a discrete mean value theorem for the Green's function. 
This latter theorem
relies on Green's function estimates in~\cite{lawler95},
and Proposition~\ref{prop-L} given in the Appendix.
The second result is Proposition~\ref{prop-uniform}, which
we prove in Section~\ref{sec-uniform}. The set we wish to
hit is not a sphere, and the proof is inspired by 
Lemma 5 of \cite{lawler92}, which only gives an upper bound.

\subsection{A discrete mean value theorem}\label{sec-gauss}
Our main result in this section is the following.
\bt{theo-gauss} Let $\{\Delta_n,n\in \N\}$ be a positive
sequence with $\Delta_n\le K n^{1/3}$ for some constant $K$, and
set $r_n=n-\Delta_n$.
There is a constant $K_a$, such that for
any $\Lambda\subset\partial\B_n$
\be{prereq-main}
\Big|E\cro{M(|\B_{r_n}|\id_0,n,\Lambda)}
-E\cro{M(\B_{r_n},n,\Lambda)}\Big| \le K_a |\Lambda|.
\ee
\et
Written explicitly, \reff{prereq-main} reads
\be{prereq1}
\big||\B_{r_n}|\times \P_0\pare{S(H_n)=z^*}
-\sum_{y\in \B_{r_n}} \P_y\pare{S(H_n)=z^*}\big|
\le K_a.
\ee
We now recall a classical decomposition 
(Lemma 6.3.6 of \cite{lawler-limic}). For a finite subset $\Lambda$,
$y\in \Lambda$, and $z^*\in \p \Lambda$
\be{har-5}
\P_{y}\left(S(H(\p \Lambda) = z^*\right)
=\frac{1}{2d}\sum_{z\in \Lambda,z\sim z^*}\!\! G_\Lambda(y,z).
\ee
By \reff{prereq1} and \reff{har-5} with $\Lambda=\B_n$, we have reduced
Theorem~\ref{theo-gauss} to proving
a discrete mean value theorem which we formulate next.
We keep the same notation as Theorem~\ref{theo-gauss}.
\bp{prop-gauss} For $z\in \B_n$, and $n-\|z\|\le 1$,
\be{gauss-estimate}
\big|\ |\B_{r_n}|\times G_n(0,z) - \sum_{y\in\B_{r_n}} G_n(y,z)\ \big|
\le K_a.
\ee
\ep
\br{rem-gauss} Note that a related (but distinct) property was also at
the heart of \cite{lawler92}. Namely, for $\epsilon>0$, and
$n$ large enough, if $z\in \B_n$, and $n-\|z\|\ge \epsilon n$,
\be{heart-lawler}
|\B_{n}|\times G_n(0,z) \ge  \sum_{y\in\B_n} G_n(y,z).
\ee
\er
\bpr
We use an improved version of Lemma 2 of \cite{lawler95}.
Using $G_n(0,z)=G(0,z)-\E_z[G(0,S(H_n))]$ (Proposition 1.5.8 of
\cite{lawler}), and \reff{main-green}, one
obtains by a Taylor expansion
that for a constant $K_1$ (independent on $n$)
\be{prereq5}
\big| \omega_d G_n(0,z)-2\frac{\alpha(z)}{n^{d-1}}\big|
\le \frac{K_1}{n^d},\quad\text{where}\quad
\alpha(z)=\E_z\cro{\|S(H_n)\|-\|z\|}.
\ee
Now, $r_n^d=n^d-d\Delta_n n^{d-1}+O(\Delta_n^2 n^{d-2})$, so that
using \reff{prereq5}, and the hypothesis $\Delta_n=O(n^{1/3})$, and
$0\le n-\|z\|\le 1$
\be{prereq6}
\begin{split}
\|\B_{r_n}\| G_n(0,z)&=\pare{r_n^d+O(r_n^{d-1})}
\pare{ 2\frac{\alpha(z)}{n^{d-1}}+O(\frac{1}{n^d})}\\
&= \pare{n^d-d\Delta_n n^{d-1}+O(\Delta_n^2 n^{d-2})+O(n^{d-1})}
\pare{2\frac{\alpha(z)}{n^{d-1}}+O(\frac{1}{n^d})}\\
&=2\alpha(z)(n-d\Delta_n)+O(1)
\end{split}
\ee
A martingale argument (Lemma 3 of \cite{lawler95}) yields
for a constant $K_l$
\be{prereq7}
\big| \sum_{y\in \B_n} G_n(y,z)-2\alpha(z)n\big|\le K_l.
\ee
Proposition~\ref{prop-L} of the Appendix reads here as follows. 
There is $K_b$ such that for $z\in \B_n$ with $n-\|z\|\le 1$
\be{prereq8}
\big| \sum_{y\in \S(r_n,n)} G_n(y,z)-2\alpha_0(z)d\Delta_n 
\big|\le  K_b,\quad\text{where}\quad
\alpha_0(z)=\E_z\cro{\|S(H_n)\|-\|z\|\big| H_n<H(B_{r_n})}.
\ee
Now, combining \reff{prereq7} and \reff{prereq8} we obtain 
(when $0\le n-\|z\|\le 1$)
\be{prereq9}
\big| \sum_{y\in \B_{r_n}} G_n(y,z)
-2n\pare{\alpha(z)-\alpha_0(z)}d\Delta_n\big|\le K_l+K_b.
\ee
Now, we combine \reff{prereq5} and \reff{prereq9} we obtain
for a constant $K_2$,
\be{prereq10}
\left| \pare{ |\B_{r_n}| G_n(0,z)-\sum_{y\in \B_{r_n}} G_n(y,z)}+
2\pare{\alpha_0(z)-\alpha(z) } d\Delta_n\right| \le K_2.
\ee
We now bound $|\alpha_0(z)-\alpha(z)|$ by the following expression
\be{prereq11}
\P_z\pare{H(B_{r_n})<H_n}\times\pare{\alpha_0(z)
+\E_z\cro{ \|S(H_n)\|-\|z\|\big| H_n>H(B_{r_n})}}.
\ee
Now, it is a classical estimate (see \reff{D1outer})
that there is $K_0$ such that for any $z\in \AA(n-1,n)$,
\be{prereq12}
\P_z\pare{H(B_{r_n})<H_n}\le \frac{K_0}{\Delta_n}.
\ee
Thus,
\be{prereq13}
\Delta_n|\alpha_0(z)-\alpha(z)|\le 
2\Delta_n \P_z\pare{H(B_{r_n})<H_n}\le 2K_0.
\ee
The desired result follows at once.
\epr
\subsection{Proof of Proposition~\ref{prop-uniform}}\label{sec-uniform}

For $j\geq 0$, consider $z_j$ in $\Sigma_j$.
We show that for all $z^*$ in ${\cal C}(z_j)$
and for suitable positive constants $\alpha_1$, $\alpha_2$, 
\be{goal1}
\frac{\alpha_1}{h_j^d}
\leq\P_{z_j}\left(S(\sigma_j) = z^*\right)
\leq\frac{\alpha_2}{h_j^d}
\ee
First, $z^*=z_j$ is a flashing position when $X_j=1$. This
happens with probability $1/h_j^d$, and gives the result.
Now, consider $z^*\in \C(z_j)\bs\{z_j\}$. We recall that
the unbiased Bernoulli $Y_j$ decides whether we flash on
$\p B(z_j,R_j)$ or on $\p \AA(r_j-R_j,r_j+R_j)$, where $R_j$
has density $g_j$ given in \reff{density}.

\noindent{\bf Step 1: Proof of the upper bound in~\reff{goal1}.}
The following obvious facts follow from Lemma~\ref{goodnn}.
\[
\text{(i)}\quad
z^*\in\p\B(z_j,\|z^*-z_j\|),\quad\text{and}\quad
z^*\in\p\AA(r_j - |\|z^*\|-r_j|, r_j +|\|z^*\|-r_j|).
\]
\[
\text{(ii)}\quad
z^*\not\in\p\B(z_j,\|z^*-z_j\|-1),\quad\text{and}\quad
z^*\not\in\p\AA(r_j - |\|z^*\|-r_j|+1, r_j +|\|z^*\|-r_j|-1).
\]
This means that if $Y_j=1$, then $R_j\in [\|z^*-z_j\|-1,\|z^*-z_j\|[$,
whereas if $Y_j=0$, then $R_j\in[|\|z^*\|-r_j|-1,|\|z^*\|-r_j|[$.
Thus, there is a constant $C$ such that
\be{har-6}
\text{(i)}\quad P\pare{Y_j=1,R_j\in [\|z^*-z_j\|-1,\|z^*-z_j\|[}
\le C\frac{\|z^*-z_j\|^{d-1}}{h_j^d},
\ee
and
\be{har-8}
\text{(ii) }\quad
P\pare{Y_j=0,R_j\in [|\|z^*\|-r_j|-1,|\|z^*\|-r_j|[}
\le C\frac{|\|z^*\|-r_j|^{d-1}}{h_j^d}.
\ee
In the case $z^*\in\p\B(z_j,\p R_j)$, the upper bound \reff{right-alex}
then follows from (i) of \reff{har-6},
and \reff{hitting-position} of Section~\ref{sec-green}.
We consider now $Y_j=0$. To simplify the notation we set for $h>0$,
\[
\D_h=\AA(r_i-h,r_i+h),\quad\text{and}\quad
\tilde \D_h=\AA(r_i-\frac{h}{2},r_i+\frac{h}{2}),
\]
and define two stopping times
\[
\tau=\inf\acc{n\ge 0:\ S(n)\in \p \D_h\cup \{z_j\}},\quad\text{and}\quad
\tau^+=\inf\acc{n\ge 1:\ S(n)\in \D_h^c\cup \{z_j\}}.
\]
It is enough to prove that for some constant $c$, and for $h$
such that $z^*\in \p \D_h$, (and $h\in  [|\|z^*\|-r_j|-1,|\|z^*\|-r_j|[$)
\be{right-alex}
\P_{z_j}\pare{ S(H(\D_h^c))=z^*}\le \frac{c}{h^{d-1}}.
\ee
We use a last exit decomposition, and the strong Markov property to get
\be{aub}
\begin{split}
\P_{z_j}\left(S(H(\D_h^c)) = z^*\right)
&\leq
G_{\D_h} (z_j, z_j) \P_{z^*}\left(H(\tilde\D_h) < \tau^+ \right)
\max_{x\in\partial\tilde\D_h^c}\P_{x}\left(S(\tau) = z_j\right)\\
& = \P_{z^*}\left( H(\tilde\D_h) < \tau^+ \right)
\max_{x\in\partial\tilde\D_h^c}G_{\D_h} (x, z_j)\\
&\le \P_{z^*}\left( H(\tilde\D_h) < \tau^+ \right)
\max_{x\in\partial\tilde\D_h^c}G(x, z_j).
\end{split}
\ee
It follows, from a Gambler's ruin estimate, that for a constant $K_0$ 
\be{flashing-12}
 \P_{z^*}\pare{H(\tilde \D_h)<\tau^+}\le \frac{K_0}{h}.
\ee
Now, from the Greens' function asymptotics \reff{main-green}
\be{flashing-13}
\sup_{x\in \p \tilde \D_h^c}
G(x,z_j)\le \sup_{x\in \p \tilde \D_h^c}\pare{
\frac{C_d}{\|x-z_j\|^{d-2}}+\frac{K_g}{\|x-z_j\|^{d}}}.
\ee
Note that the distance between $z_j$ and $\tilde \D_h$ is of order $h$.
We use \reff{flashing-12} and
\reff{flashing-13} in \reff{aub} to obtain \reff{right-alex}.

\noindent{\bf Step 2: Proof of the lower bound in~\reff{goal1}.}
Note the following two facts.
\begin{itemize}
\item[(iii)] 
By Lemma~\ref{goodnn}, $z^*$ has a nearest neighbor, say $z$, in
$\B(z_j,h)$ with 
\[
\|z-z_j\|\le \|z^*-z_j\|-\frac{1}{4\sqrt{d}}.
\] 
This means that if
$h\in [\|z^*-z_j\|-1/(4\sqrt{d}),\|z^*-z_j\|[$, then
$z^*\in \p\B(z_j,h)$.
\item[(iv)] 
By Lemma~\ref{goodnn}, $z^*$ has a nearest neighbor, say $z$, in
$\AA(r_j-h,r_j+h)$ with 
\[
|\|z\|-r_j|\le |\|z^*\|-r_j|-\frac{1}{4\sqrt{d}}.
\]
This means that if
$h\in [|\|z^*\|-r_j|-1/(4\sqrt{d}),\|z^*\|-r_j|[$, then
$z^*\in \p\D_h$.
\end{itemize}

We deal separately with the cases $|\|z^*\| - r_j| < h_j /2$
and $|\|z^*\| - r_j| \geq h_j /2$.

Consider first the case $|\|z^*\| - r_j| < h_j /2$. On the event $Y_j=1$,
and for $h$ such that $z^*\in \p \B(z_j,h)$, we have
\be{har-7}
P(Y_j=1,\ R_j\in [\|z^*-z_j\|-1/(4\sqrt{d}),\|z^*-z_j\|[)\ge 
\frac{ch^{d-1}}{h_j^d}
\P_{z_j}\pare{ S(H(\p \B(z_j,h)))=z^*}\ge \frac{c}{h_j^{d-1}}.
\ee
Thus, for some constant $\alpha_1$ (that depends on $d$), we have
\reff{goal1}.

Consider now the case $|\|z^*\| - r_j| \geq h_j /2$. 
It is enough to prove, for $h$ such that $z^*\in \p\D_h$,
and for some constant $c$ (that depends on $d$)
\be{goal3}
P_{z_j}\left(S(H(\D_h^c))= z^*\right)
\geq \frac{c}{h_j^{d-1}}.
\ee
Let $y^*$ be the closest site of $\p \B(0,r_j)$ to the segment $[0,z^*]$,
and $x^*$ be the closest site of $\p \B(0,r_j+h/2)$
to the segment $[0,z^*]$.
We set $\Gamma = \B(x^*, \|z^* - x^*\|)\cap \tilde \D_h$. 
It may be that 
$\B(x^*, \|z^* - x^*\|)\cap \D_h^c=\not\emptyset$,
and if so, one would only have to consider a site at a distance
1 from $z^*$, say $\tilde z\in \D_h$, and such that
$\B(x^*, \|\tilde z - x^*\|)\cap \D_h^c=\emptyset$, and work with
$\tilde z$ instead of $z^*$ in the sequel.
We assume henceforth
that $\B(x^*, \|z^* - x^*\|)\cap \D_h^c=\emptyset$. 

By~\reff{har-5} with $\Lambda=\D_h$, and the strong Markov property,
\be{har-1}
\begin{split}
\P_{z_j}\left(S(H(\p\D_h)) = z^*\right)&\ge
G_{\D_h} (z_j, z_j) \P_{z^*}\left(H(\Gamma) < \tau^+ \right)
\min_{x\in\Gamma}\P_{x}\left(S(\tau) = z_j\right)\\
&\ge \P_{z^*}\left(H(\Gamma) < \tau^+\right)
\min_{x\in\Gamma}G_{\D_h} (x, z_j).
\end{split}
\ee
Since $z^*\in{\cal C}(z_j)$, $y^*$ and $z_j$
can be connected by 20 overlapping balls of radius $h_j / 10$
in such a way that, applying Harnack's inequality 20 times
(see Theorem 6.3.9 in~\cite{lawler-limic})
to the harmonic functions $G_{\D_h}(x,\cdot)$,
we can estimate from below the last factor
in \reff{har-1}. There is a constant $K_2$ such that
\be{har-2}
\begin{split}
\min_{x\in\Gamma\cap\partial_-\D_h'}G_{\D_h} (x, z_j)
\ge& c_H^{20}
\min_{x\in\Gamma\cap\partial_-\D_h'}G_{\D_h} (x, y^*)\\
\ge &  
c_H^{20}\min_{x\in\Gamma\cap\partial_-\D_h'}G(x, y^*)
- \E_x\left[G(S(H(\D_h^c)), y^*)\right]\\
\ge&
\frac{c_H^{20} C_d}{h^{d-2}}\pare{\frac{1}{(\sqrt{2}/2)^{d-2}}-1}
\ge\frac{K_2}{h_j^{d-2}}.
\end{split}
\ee
As a consequence of \reff{har-2}, we just need to prove
that the first factor in~\reff{har-1} is of order $1/h_j$ at least.
We realize the event $\{ H(\Gamma)<\tau^+\}$ in two moves: we first
hit the sphere $\B(x^*,R^*/2)$, and then we exit from the cap
$\p \B(x^*,R^*)$ which lies in $\tilde \D_h$.
\be{flashing-18}
\begin{split}
\P_{z^*}\pare{ H(\Gamma)<\tau^+}\ge &\P_{z^*}
\pare{ H(\B(x^*,R^*/2))<H(\B^c(x^*,R^*))}\\
&\inf_{y\in \p \B(x^*,R^*/2)} \P_y\pare{
H\pare{\B^c(x^*,R^*)\cap \tilde \D_h}= H\pare{\B^c(x^*,R^*)}}
\end{split}
\ee
We invoke again Harnack's inequality to have for $y\in \p B(x^*,R^*/2)$
\be{flashing-19}
\P_y\pare{
H\pare{ \B^c(x^*,R^*)\cap \tilde \D_h}= H\pare{\B^c(x^*,R^*)}}\ge c_H
\P_0\pare{
H\pare{E\cap \B^c(x^*,R^*)}= H\pare{\B^c(x^*,R^*)}}.
\ee
We invoke now \reff{hitting-position} to obtain for some constant $K_3$
\be{flashing-20}
\P_0\pare{H\pare{\B^c(x^*,R^*)\cap \tilde \D_h}= H\pare{\B^c(x^*,R^*)}}
\ge c_1\frac{| \p\B(x^*,R^*)\cap \tilde \D_h|}{|\p\B(x^*,R^*)|}\ge K_3.
\ee
We gather now \reff{flashing-18},
\reff{flashing-19} and \reff{flashing-20} to obtain the desired
lower bound.
\section{The flashing process fluctuations}\label{sec-results}

In this section we prove Theorems~\ref{theo-flashing} 
and~\ref{theo-optimal}.
To do so we use the construction in terms of exploration waves
of Section~\ref{sec-wave}.

\subsection{Tiles}
We recall that we have defined a {\it cell} of $\S_j$
in \reff{def-cell}, as the intersection of a cone with $\S_j$.
We need also a smaller structure.
We define, for any $z_j$ in $\Sigma_j$,
\begin{equation}
\tilde{\cal C}(z_j) = {\cal S}_j \cap\left\{
  x \in {\mathbb R}^d :\:  
  \exists \lambda \geq 0, \exists y\in B(z_j, h_j/5), x = \lambda y
\right\}.
\label{dfn-cell2}
\end{equation} 
As in Lemma 12 in~\cite{lawler95}, concerning locally finite coverings,
we claim that, for $h_o$ large enough,
there exist a positive constants $c_1$,
and, for each $j\geq 0$, a subset $\tilde\Sigma_j$
of $\Sigma_j$ such that 
\begin{equation}
\left| \tilde\Sigma_j \right|
\leq c_1\frac{|\Sigma_j|}{h_j^{d-1}}
\quad\mbox{and}\quad
S_j = \bigcup_{z_j\in\tilde\Sigma_j} \tilde{\cal C}(z_j).
\label{local-covering}
\end{equation}
For any $z_j\in \Sigma_j$,
we call {\em tile centered at }$z_j$, 
the intersections of $\tilde{\cal C}(z_j)$ with $\Sigma_j$.
We denote by $\T(z_j)$  a {\em tile centered at }$z_j$, and
 by ${\cal T}_j$
the set of tiles associated with the shell ${\cal S}_j$:
\begin{equation}
{\cal T}_j = \left\{ {\cal T}(z_j):\:\ z_j\in\tilde\Sigma_j \right\}.
\end{equation}
Let us explain the reason for $h_j/5$ in the definition of a tile.
It implies a fundamental feature of the flashing process. For any
$z\in \S_j$, there is $\tilde z_j\in \tilde \Sigma_j$ such that
\be{feature-coupon}
z\in \bigcap \acc{\C(y): \ y\in \T(\tilde z_j)}.
\ee
Indeed, let $z_j\in \Sigma_j$ be the site realizing
the minimum of $\{\|z-y\|:\ y\in \Sigma_j\}$. There is
$\lambda>0$ and $u\in B(z_j,1)$, such that $z=\lambda u$.
Now, there is $\tilde z_j\in \tilde \Sigma_j$ such that
$\|\tilde z_j-z_j\|<h_j/5$, and for any $y\in \T(\tilde z_j)$,
we have $\|y-z_j\|<2h_j/5$. Thus
\[
u\in \bigcap \acc{B(y,\frac{h_j}{2}): \ y\in \T(\tilde z_j)}\Longrightarrow
z\in \bigcap \acc{\C(y): \ y\in \T(\tilde z_j)}.
\] 
\subsection{The inner ball}
For $n\geq 0$, we take $N = |{\mathbb B}_n|$, we recall that
${\cal A}^*(N) = \cup_{k\geq 1} {\cal A}^*_k(N)$,
and write ${\cal A}^*$ instead of ${\cal A}^*(N)$. We consider
\begin{equation}
T^* = \min\left\{
  k\geq 1 :\: \cup_{j<k}{\cal S}_j \not\subset {\cal A}^*_k
\right\}.
\end{equation}
We have, for $l$ with $r_l <n$,
\begin{equation}
P\left({\mathbb B}(0,r_l + h_l)\not\subset{\cal A}^*\right)
\leq
\sum_{k\leq l} P(T^* = k+1)
\label{est0}
\end{equation}
and we estimate from above the probability $P(T^* = k+1)$
assuming $r_k < n$.
For $k\geq 1$ and $\Lambda\subset\Sigma_k$,
we call $W_k(\Lambda)$ the number of unsettled explorers
that stand in $\Lambda$ after the $k$-th wave,
that is
\begin{equation}
W_k(\Lambda)
= \sum_{i = 1}^N \id_{\Lambda}\left(\xi_k(i)\right).
\end{equation}
We now look at the
{\it crossings} of tiles of $\T_k$. On the one hand, we will use that
if $W_k({\cal T})$ is {\it large}, then it is unlikely that a hole
appears in the cell containing $\T$. On the other hand, if $r_k$ is
{\it small} it is unlikely that $W_k({\cal T})$ is {\it small}.
We make precise what we intend by {\it small} and {\it large}.
For this purpose, we will show in \reff{condtalt}
that for some constant $\kappa_1>0$, and any tile $\T\in \T_k$
\be{def-h}
E[W(\T)]\ge \kappa_1 (n-r_k) h_k^{d-1},
\quad\text{and we define}\quad h = n^{\frac{1}{d+1}}\ge
\sup_{k:r_k\le n} h_k.
\ee
For any positive constant $A$, we write
\be{est1}
\begin{split}
P\left(T^* = k+1\right)
=& P\left( T^* = k + 1,
  \exists {\cal T}\in{\cal T}_k, W_k({\cal T}) < Ah_k^d\log n
\right)\nonumber\\
&+ P\left( T^* = k + 1,
  \forall {\cal T}\in{\cal T}_k, W_k({\cal T}) \geq Ah_k^d\log n
\right),
\end{split}
\ee
and we estimate separately each term in the right hand side of \reff{est1}.

\paragraph{Estimating the first term.}
We show here that (for $\kappa_1$ and $h$ appearing in \reff{def-h}) for
$k$ such that
\be{inter-10}
r_k\leq n - \frac{2A}{\kappa_1} h\log n,
\ee
there is a constant $\kappa_2>0$, and $n$ large enough, such that
\be{inner-1}
P\left( T^* = k + 1,
\exists {\cal T}\in{\cal T}_k, W_k({\cal T}) < Ah_k^{d}\log n \right)
\le |\S_k| \exp\pare{-\kappa_2 A^2 \log^2n}.
\ee
On $\{T^* = k + 1\}$, we have $A_k^*=\B(0, r_k - h_k)$.
On $\{T^* = k + 1\}$, and for any $\T\subset \T_k$, we consider
a variable $L_k(\T)=M(\B(0, r_k - h_k), r_k, \T)$ independent of
$W_k(\T)$, and define $M_k(\T)=W_k(\T) + L_k( \T)$.
We have the equality in law, on $\{T^* = k + 1\}$,
\begin{equation}\label{identity}
M_k(\T)\stackrel{\text{law}}{=}M(N\id_{\{0\}}, r_k, \T),\quad\text{and}
\quad W_k({\cal T}) = M_k({\cal T}) - L_k({\cal T}).
\end{equation}
As a consequence
\begin{equation}
P\left(
  T^* = k + 1,
  \exists {\cal T}\in{\cal T}_k, W_k({\cal T}) < Ah_k^d\log n
\right)
\leq |\tilde\Sigma_k| \max_{{\cal T}\in {\cal T}_k}
P\left(
  M_k({\cal T}) - L_k({\cal T})< Ah_k^d\log n
\right).
\label{maj0}
\end{equation}
\par\noindent
We first estimate the number of explorers stopped on
a tile $\cal T$ of ${\cal T}_k$.

Now, $M_k({\cal T})$ and $L_k({\cal T})$ are dependent
random variables, but both are sums of independent Bernoulli
variables, for which Lemma~\ref{lem-BL} is designed.
We introduce two notations.
For a variable $X$, let $\bar X=X-E[X]$, and let
\be{def-xk}
2 \bar x_k=E\cro{M_k({\cal T}) - L_k({\cal T})}-Ah_k^d \log(n).
\ee
Since we need $\bar x_k$ of \reff{def-xk} to be positive,
we will choose $A$ and $k$ such that
\begin{equation}
\label{condt}
E\cro{M_k({\cal T}) - L_k({\cal T})}\ge 2 Ah_k^d \log(n),
\quad\text{which implies}\quad \bar x_k\ge \frac{1}{2} Ah_k^d \log(n).
\end{equation}
Then,
\be{inter-1}
P\pare{M_k({\cal T}) - L_k({\cal T})< Ah_k^d\log n} =
P\pare{\bar M_k({\cal T}) - \bar L_k({\cal T})< -2\bar x_k}.
\ee
In order to estimate $E[M_k({\cal T}) - L_k({\cal T})]$, we
invoke Theorem~\ref{theo-gauss}, with $n=r_k$, and $\Delta_n=h_k$
(the hypothesis $h_k=O(r_k^{1/3})$ holds here). We have
for some positive constants $\kappa'$, $\kappa_1$, and for $n$ large enough
\be{condtalt}
\begin{split}
E\left[ M_k({\cal T}) - L_k({\cal T}) \right] = &
E[M((|{\mathbb B}_n|-|{\mathbb B}_{r_k-h_k}|)\id_{0},r_k,\T)]\\
&\quad + E[M(|{\mathbb B}_{r_k-h_k}|\id_{0},r_k,\T)]-
E[M({\mathbb B}_{r_k-h_k},r_k,\T)]\\
\ge & \left( |{\mathbb B}_n| - |{\mathbb B}_{r_k-h_k}|
\right) {\mathbb P}_0\left( S(H_k)\in{\T} \right) -O(h_k^{d-1})\\
\geq&
\kappa'(n^d - (r_k-h_k)^d)\frac{h_k^{d-1}}{r_k^{d-1}}-O(h_k^{d-1})\\
\geq& \kappa_1(n-r_k)h_k^{d-1}.\qquad(\text{recall that }\ n-r_k>h).
\end{split}
\ee
Note that the ultimate inequality in \reff{condtalt} is
the estimate in \reff{def-h}.
In view of \reff{condtalt}, condition~(\ref{condt}) is ensured
if $A$ and $k$ satisfy \reff{inter-10}. 

Note that \reff{condtalt} implies that $E[M_k({\cal T})]
\ge E[L_k({\cal T})]$, so that \reff{eq-BL2} of Remark~\ref{rem-BL}
requires only an upper bound on $E[M_k({\cal T})]$. Thus, we only
treat the latter quantity.

We distinguish two cases: (i) when $r_k$ is {\it close} to $n$, 
(ii) when $r_k$ is {\it small} compared to $n$.

\noindent{\bf Step 1:} We assume $n-h\ge r_k \geq n /2$.\\
We set here $2x_k=Ah_k^d\log n$. \reff{condt} and \reff{inter-1} 
imply that
\be{inter-11}
P\pare{M_k({\cal T}) - L_k({\cal T})< Ah_k^d\log n}
\le P\pare{\bar M_k({\cal T})<-x_k} +P\pare{\bar L_k({\cal T})>x_k}.
\ee
To be in the CLT regime of Lemma~\ref{lem-BL} when dealing
with the right hand side of \reff{inter-11}, we need
\be{inter-2}
0<x_k< E[M_k({\T})].
\ee
Let us now estimate $E[M_k(\T)]$.
Note that by using \reff{hitting-position} and $r_k\ge n/2$,
we have for positive constants $K_1,K_1'$
\be{estvar}
\begin{split}
K_1' n^d\pare{\frac{h_k}{r_k}}^{d-1}\ge&
E[M_k(\T)]=|\B_n| \P_0(S(H_k)\in \T)\\
\ge & K_1 n^d\pare{\frac{h_k}{r_k}}^{d-1}
\geq K_2 n^d\frac{h_k^{d-1}}{n^{d-1}} = K_2 nh_k^{d-1}.
\end{split}
\ee
Thus, $E[M_k({\cal T})]> x_k$ and
we are in the Gaussian regime for \reff{eq-BL}.
Similarly, as in \reff{estvar}, we have a lower bound
\be{inter-5}
E[M_k(\T)]\le K_2' nh_k^{d-1}.
\ee
Thus, there is $\kappa_2>0$ such that for a large enough $n$,
\be{inter-4}
P\pare{\bar M_k({\cal T}) \leq -x_k}\le e^{-\frac{x_k^2}{4 E[M_k(\T)]}}
\le \exp\pare{-\frac{A^2h_k^{2d}\log^2n}{16 K_2'nh_k^{d-1}}}
=\exp\pare{-\kappa_2 A^2 \log^2n}.
\ee
As already noted, \reff{eq-BL2} yields
\be{bound-L}
P(\bar L_k(\T)\ge x_k)\le\exp(-\kappa_2 A^2 \log^2n).
\ee
\par\noindent
\noindent{\bf Step 2:} We assume $r_k < n/2$.\\
We have, using \reff{condtalt}, for $n$ large enough
\be{inter-7}
2\bar x_k\geq \kappa'(n^d-r_k^d) \pare{\frac{h_k}{r_k}}^{d-1}-Ah_k^d\log n
\geq \frac{\kappa'}{2}n^d \pare{\frac{h_k}{r_k}}^{d-1}
-Ah_k^d\log n \ge \frac{\kappa'}{4}n^d \pare{\frac{h_k}{r_k}}^{d-1}.
\ee
We define here
\be{def-xkbis}
x_k=\frac{\kappa'}{16}n^d \pare{\frac{h_k}{r_k}}^{d-1},\quad
(\text{and note that }\quad x_k\ge \frac{\kappa'}{16}n h_k^{d-1}).
\ee
As previously, we have \reff{inter-11}.
From \reff{estvar}, we have for some positive $K_2'$
\be{inter-8}
E[M_k({\cal T})]\leq K_2' n^d \pare{\frac{h_k}{r_k}}^{d-1}.
\ee
Now, using Lemma~\ref{lem-BL}, for $n$ large enough
\be{maj4}
\begin{split}
P\left(M_k({\cal T}) - L_k({\cal T})< Ah_k^d\log n \right)
&\leq 2\exp\pare{-\frac{x_k}{4} \min\pare{1,\frac{\kappa'}{8K_2'}}}\\
&\le \exp\pare{-\kappa_2 A^2 \log^2n}.
\end{split}
\ee
\par\noindent
Collecting~(\ref{inter-4}), \reff{bound-L} and~(\ref{maj4})
together with~(\ref{maj0}) and~(\ref{local-covering}),
we conclude that \reff{inner-1} holds.

\paragraph{Estimating the last term.}
The last term in the right hand side of \reff{est1}
is bounded using a simple coupon-collector argument.
Indeed, the event $\{T^* = k + 1\}$ means that there is one
uncovered site in $\S_k$. By \reff{feature-coupon}, there is $z_k\in
\tilde \Sigma_k$, such that this site is a possible settling
position of all explorers stopped in $\T(z_k)$. Now, if
$\{W_k(\T(z_k)) \geq Ah_k^d\log n\}$, Proposition~\ref{prop-uniform}
tells us that the probability of not covering this site is
less than $(1-\alpha_2/h_k^d)$ to the power $Ah_k^d \log(n)$.
In other words,
\be{inner-2}
\begin{split}
P\left( T^* = k + 1,
  \forall {\cal T}\in{\cal T}_k, W_k({\cal T}) \geq Ah_k^d\log n \right)
&\leq
|{\cal S}_k|\left(
  1 - \frac{\alpha_2}{h_k^d} \right)^{Ah_k^d\log n}\\
&\leq |{\cal S}_k| \exp\left\{ -\alpha_2 A\log n \right\}.
\end{split}
\ee

\paragraph{Conclusion.}
First, choose $A$ large enough so that
\be{inter-9}
|{\mathbb B_n}| \exp\pare{ -\alpha_2 A\log n }\le \frac{1}{n^2}.
\ee
Recall the decomposition~(\ref{est0} and~(\ref{est1}), and assume
that $r_l$ satisfies \reff{inter-10}.
Then, \reff{inner-1} and \reff{inner-2} yield that for $n$ large enough
\begin{eqnarray}
P({\mathbb B}(0,r_l + h_l \not\subset {\cal A}^*(|{\mathbb B_n}|))
\leq
|{\mathbb B_n}| \left(\exp\pare{-\kappa_2 A^2\log^2 n} +
 \exp\pare{ -\alpha_2 A\log n } \right)
\le \frac{2}{n^2}.  \label{cl}
\end{eqnarray}
The right-hand side in~(\ref{cl})
is summable, and Borel-Cantelli lemma yields
the inner control of Theorem~\ref{theo-flashing}.

\subsection{The outer ball}\label{sec-outer}

This section follows closely~\cite{lawler95}.
The features of the flashing
process allow for some simplification.
We keep the notation of the previous section.
There, we proved that for some positive constant $\delta$
\[
P\pare{
\left\{{\mathbb B}(0, n - \delta h\log n)\subset {\cal A}^* \right\}}
=1 - \epsilon^\delta(n),\quad\text{with}\quad
\sum_{n\geq 1} \epsilon^\delta(n) < +\infty.
\]
As consequence, the following conditional law
can be seen as a slight modification of $P$.
\begin{equation}
P^\delta(\cdot) = P\pare{\cdot |
\left\{{\mathbb B}(0, n - \delta h\log n)\subset {\cal A}^* \right\}}.
\end{equation}

We begin by proving that, under $P^\delta$,
the probability to find some $k$
with $r_k < 2n$ and some tile ${\cal T}$
in ${\cal T}_k$ with $W_k({\cal T})$
larger than or equal to $2Ah^d\log n$
for a large enough $A$ decreases
faster than any power of $n$.
First, note that, under $P^\delta$, we have
\begin{equation}
W_k({\cal T}) \leq M_k({\cal T}) - L_k^\delta({\cal T}),\quad\text{with}
\quad 
L_k^\delta = M\left(
  {\mathbb B}(0, n-\delta h \log n), r_k, {\cal T}
\right).
\end{equation}
Now,
\begin{equation}
P^\delta\left(W_k({\cal T}) \geq 2Ah^d\log n\right)
\leq
P^\delta\left(M_k({\cal T}) - L_k^\delta({\cal T})\geq 2Ah^d\log n\right).
\end{equation}
By Theorem~\ref{theo-gauss},
for some positive constants $K'$, $K$ and for $n$ large enough
\be{inner-3}
\begin{split}
E\left[ M_k({\cal T}) - L_k^\delta({\cal T})
\right]
&\leq
K'\left(n^d - (n-\delta h\log n)^d\right)\frac{h_k^{d-1}}{r_k^{d-1}}
+ O(h_k^{d-1})\\
&\leq
K'n^{d}d\frac{\delta h\log n}{n}\frac{h_k^{d-1}}{r_k^{d-1}}
+ O(h_k^{d-1})\\
&\leq
K'd\delta hh_k^{d-1}\log n
+ O(h_k^{d-1}) \leq Kh^{d}\log n.
\end{split}
\ee
Choosing $A\geq K$, we get for $n$ large enough so that
$P({\mathbb B}(0, n - \delta h\log n)\subset {\cal A}^*)\ge 1/2$
\begin{eqnarray}
P^\delta\left(W_k({\cal T}) \geq 2Ah^d\log n\right)
&\leq& 2 P\left(
  M_k({\cal T}) - E[M_k({\cal T})] 
  \geq \frac{A}{2}h^d\log n
\right)\nonumber\\
&& + 2P\left(
  L_k^\delta({\cal T}) - E[L_k^\delta({\cal T})] 
  \leq -\frac{A}{2}h^d\log n.
\right)
\end{eqnarray}
As in the previous section 
$E[M_k({\cal T})]$ 
is of order $n^dh_k^{d-1}/r_k^{d-1}$, i.e,
of order $nh^{d-1}$.
In addition $E[L_k^\delta({\cal T})]$
is smaller than $E[M_k({\cal T})]$.  
We conclude once again by invoking Lemma~\ref{lem-BL}.

Now, let $F_k$ denote the event
that no tile ${\cal  T}$ in $\Sigma_k$ contains
more than $2Ah^d\log n$ unsettled explorers
after the $k$-th  exploration wave. 
We denote by $\FF_k=\sigma(\xi_0,\dots,\xi_k)$, and note
that $F_k$ and
$\{\B(0,n-\delta h^d \log(n))\subset \A^*\}$ are $\FF_k$-measurable.

For any tile
${\cal T}\in \T_k$, let $z_k\in \tilde \Sigma_k$ be such that
$\T=\T(z_k)$ and denote by $\tilde {\cal C} = \tilde{\cal C}(z_k)$.
We are entitled, by Proposition~\ref{prop-uniform},
to use a coupon-collector estimate on the number of
settled explorers during the $k+1$-th exploration wave.
On $F_k\cap \{\B(0,n-\delta h^d \log(n))\subset \A^*\}$, and
for some positive constant $K_1$
\be{inner-4}
\begin{split}
E\left[ {\cal A}^*_{k+1}\cap \tilde{\cal C}
  \Big| \FF_k \right]
&\geq
|\tilde{\cal C}|\left( 1 - \left( 1 - \frac{\alpha_1}{h_k^d}
  \right)^{W_k({\cal T})} \right)\\
&\geq
|\tilde{\cal C}|\left( 1 - \exp\left\{
    -\alpha_1\frac{W_k({\cal T})}{h_k^d} \right\}\right)\\
&\geq
\frac{|\tilde{\cal C}|}{h_k^d} W_k({\cal T})
\frac{h_k^d}{W_k({\cal T})}\left(
  1 - \exp\left\{
    -\alpha_1\frac{W_k({\cal T})}{h_k^d}
  \right\}
\right)\\
&\geq K_1 W_k({\cal T}) \inf_{x\leq 2A\log n} \frac{1 - e^{-\alpha_1x}}{x}.
\end{split}
\ee
We now write for some positive constant $K_2$
\[
\begin{split}
\inf_{x\leq 2A\log n} \frac{1 - e^{-\alpha_1x}}{x}
&\geq \frac{1}{2A\log n}\; \inf_{x\leq 2A\log n}
\frac{1 - e^{-\alpha_1x/2A\log n}}{x/2A\log n}\\
&\geq \frac{1}{2A\log n}\; \inf_{x\leq 1} \frac{1 - e^{-\alpha_1x}}{x}
\ge\frac{K_2}{\log n}.
\end{split}
\]
We conclude that on $F_k\cap \{\B(0,n-\delta h^d \log(n))\subset \A^*\}$,
\begin{eqnarray}
E\left[
  {\cal A}^*_{k+1}\cap \tilde{\cal C} \Big| \FF_k \right]
&\geq&
K_1K_2 \frac{W_k({\cal T})}{\log n} .
\end{eqnarray} 
Summing over all tiles we get, for a different constant $K$,
(because of the finite, $k$-independent overlapping between tiles),
we obtain on $F_k\cap \{\B(0,n-\delta h^d \log(n))\subset \A^*\}$
\[
E\left[ {\cal A}^*_{k+1}\cap{\cal S}_k \Big| \FF_k \right]
\geq
K\frac{W_k({\cal S}_k)}{\log n}.
\]
Also, since $W_k({\cal S}_k) \leq |{\mathbb B}(0,n)|$,  
\[
E\left[\id_{F_k\cap \{\B(0,n-\delta h^d \log(n))\subset \A^*\}}
 {\cal A}^*_{k+1}\cap{\cal S}_k\right]
\geq
K\frac{E[\id_{\{\B(0,n-\delta h^d \log(n))\subset \A^*\}}
W_k({\cal S}_k)]}{\log n} - n^dP(F_k^c).
\]
Since $P(\{\B(0,n-\delta h^d \log(n))\subset \A^*\})\ge 1/2$,
\begin{eqnarray}
E^\delta\left[ {\cal A}^*_{k+1}\cap{\cal S}_k
\right]
&\geq&
K\frac{E^\delta[W_k({\cal S}_k)]}{\log n}
- 2 n^dP(F_k^c).
\end{eqnarray}  
In other words, noting that $ {\cal A}^*_{k+1}\cap{\cal S}_k=
W_k({\cal S}_k)-W_{k+1}({\cal S}_{k+1})$
\be{inner-7}
E^\delta\left[
  W_{k+1}({\cal S}_{k+1})
\right]
\leq \left(1 - \frac{K}{\log n}\right) E^\delta\left[
  W_{k}({\cal S}_{k}) \right] + 2n^dP(F_k^c).
\ee
By iterating \reff{inner-7}, we obtain
that for any $\epsilon$, $E^\delta[W_{l+
\epsilon\log^2n}({\cal S}_{l+\epsilon\log^2 n})]$,
decreases faster than any power of $n$,
when $l$ the lowest index for which $r_l \geq n$.
Also, the probability (under $P$ !) of seeing
at least one explorer reaching the shell ${\cal S}_{l+\epsilon\log^2 n}$
is summable.
Using Borel-Cantelli lemma,
this yields the proof of Theorem~\ref{theo-flashing}.

\subsection{Optimality of the fluctuation exponent}\label{sec-optimal}

Let time $k$ be such that $r_k=n-Ah$, for a large arbitrary constant $A$.
We show that $P(T^*=k+1)$ decays faster than any polynomial in $n$.

On the event $\{T^*=k+1\}$, we have, after the $k$-th wave
and for some constant $K$,
\be{optimal1}
|\B_{r_k}|=M(\B_{r_k}, r_k, \Sigma_k)=
M(A_k^*, r_k, \Sigma_k)\Longrightarrow W_k(\Sigma_k)=|\B_n|-|\B_{r_k}|\le AK
n^{d-1}\times h.
\ee
This means that there exists $z_k\in \Sigma_k$
such that, for some positive constant $K'$, 
\be{optimal2}
W_k(B(z_k, 3h)\cap \Sigma_k)
\le K'n^{d-1}\times h\times
\frac{h^{d-1}}{n^{d-1}}\le K'h^d.
\ee
By construction, only the explorers stopped
inside $B(z_k, 3h)$ can cover $\tilde{\cal C}(z_k)$ (for $n$ large enough).
As a consequence, we can think of a coupon-collector problem, where
an album of size $|\tilde{\cal C}(z_k)|$ has to be filled 
when we collect no more than $K'h^d$ coupons. 
Thus, the probability
of $\{T^*=k+1\}$ is bounded from above by the probability of filling
such an album, which is less than $\exp(-c(A) h_n^{d/2})$, for some
explicit constant $c(A)$ dependent on $A$. 
This result is based on the following simple
coupon-collector lemma (together with Proposition~\ref{prop-uniform}), 
which we did not find in the vast literature on such problems.

\bl{lem-coupon}
Consider an album of $L$ items
for which are bought independent
random cou\-pons, each of them covering one (or possibly none)
of the possible $L$ items.
If $Y_i$ is the item 
associated with the $i$-th coupons, 
we assume that for positive constants $\alpha_1,\alpha_2$, such that
for any $j=1,\dots,L$,
\be{coupon1}
\frac{\alpha_1}{L} \le P(Y_i=j)\le \frac{\alpha_2}{L}.
\ee
Let $\tau_L$ be the number of coupons needed to complete the album.
Then, for any $A>0$,
\be{coupon-main}
P(\tau_L<A L)\le 
\exp\left(-\frac{\alpha_1^2 A^2e^{-2\alpha_2 A}}{4}\sqrt{ L}\right).
\ee
\el
\bpr
We denote by $\sigma_i$ the time needed to collect the $i$-th
distinct item after having collected $i-1$
distinct items. 
The sequence $\{\sigma_1,\sigma_2,\dots,\sigma_L\}$
is not independent, but if $\GG_k=\sigma(\{Y_1,\dots,Y_k\})$,
and $\tau(k)=\sigma_1+\dots+\sigma_k$, then for $i=1,\dots,L$
\be{coupon2}
\left(1-\frac{\alpha_1(L-i+1)}{L}\right)^k\ge
P(\sigma_{i}>k\big| \GG_{\tau(i-1)})\ge 
\left(1-\frac{\alpha_2(L-i+1)}{L}\right)^k.
\ee
Indeed, calling ${\cal E}(i-1)$ the set of the first $i-1$ collected items,
\be{coupon3}
\begin{split}
P(\sigma_{i}>k\big| \GG_{\tau(i-1)})=&P\left(
\{Y_{\tau(i-1)+1},\dots, Y_{\tau(i-1)+k}\}\subset {\cal E}(i-1)
| \GG_{\tau(i-1)}\right)\\
=&\left(P\left( Y \in{\cal E}(i-1)| \GG_{\tau(i-1)}\right)\right)^k\\
= &\left( 1-P\left( Y \not\in{\cal E}(i-1)| \GG_{\tau(i-1)}\right)\right)^k.
\end{split}
\ee
Using \reff{coupon1} we deduce \reff{coupon2} from \reff{coupon3}.
Now, \reff{coupon2} gives that
\be{coupon4}
 \frac{L}{\alpha_1(L-i+1)}\ge
E[\sigma_{i}| \GG_{\tau(i-1)}]\ge \frac{L}{\alpha_2(L-i+1)},
\ee
as well as
\be{coupon11}
E[\sigma_{i}^2| \GG_{\tau(i-1)}]\le 2\frac{L^2}{\alpha_1^2(L-i+1)^2}.
\ee
Now, we look for $B$ such that
\be{coupon5}
\sum_{i={\sqrt L}}^{B{\sqrt L}} E[\sigma_{L-i}]\ge 2 AL.
\ee
Note that
\[
\sum_{i={\sqrt L}}^{B{\sqrt L}} E[\sigma_{L-i}]\ge
\frac{L}{\alpha_2}\sum_{i={\sqrt L}}^{B{\sqrt L}}\frac{1}{i+1}\ge
 \frac{L}{\alpha_2} \log(B).
\]
Thus, condition \reff{coupon5} holds for $B\ge \exp(2\alpha_2 A)$.
Finally, note that
\[
\max\left\{ E[\sigma_{L-i}| \GG_{\tau(L-i-1)}],
\ i={\sqrt L},\dots,B{\sqrt L}\right\}\le
\frac{\sqrt L}{\alpha_1},
\]
and set
\[
X_i=\frac{E[\sigma_{L-i}| \GG_{\tau(L-i-1)}]
-\sigma_{L-i}}{\left( \sqrt L/\alpha_1\right)}\le 1.
\]
For $x\le 1$, note that $e^x\le 1+x+x^2$ to obtain for $0\le \lambda\le
1$, by successive conditioning
\be{coupon7}
\begin{split}
P\left(\sum_{i={\sqrt L}}^{B{\sqrt L}} \sigma_{L-i}\le AL\right)\le &
P\left( \sum_{i={\sqrt L}}^{B{\sqrt L}} X_i\ge \alpha_1 A \sqrt L\right)\\
\le & e^{-\lambda \alpha_1 A\sqrt{ L}}\prod_{i={\sqrt L}}^{B{\sqrt L}}
\left(1+\lambda^2 \sup E[X_i^2| \GG_{\tau(L-i-1)}]\right)\\
\le & \exp\left(-\lambda \alpha_1 A\sqrt L+\lambda^2
\sum_i \sup E[X_i^2| \GG_{\tau(L-i-1)}]\right).
\end{split}
\ee
Finally,
we have, using~(\ref{coupon11}),
\be{coupon8}
\sum_{i={\sqrt L}}^{B{\sqrt L}} \sup E[X_i^2| \GG_{\tau(L-i-1)}]\le
\sum_{i={\sqrt L}}^{B{\sqrt L}} \alpha_1^2\sup 
\frac{E[\sigma^2_{L-i}|\GG_{\tau(L-i-1)}]}{L}
\le 2 B \sqrt L.
\ee
The results follows as we optimize on $\lambda\le 1$ in the upper bound
in \reff{coupon7}.
\epr
\appendix

\section{Time spent in an annulus (By S.Blach\`ere)}\label{sec-blachere}

This appendix is devoted to an asymptotic expansion of the expected 
time spent in an annulus $\A(r_n,n)$ for $r_n<n$, 
when the random walk is started at some point $z$ within the annulus, 
and before it exits the outer shell.

\bp{prop-L}
Consider a sequence $\{\Delta_n,n\in \N\}$ with
$\Delta_n\le K n^{1/3}$ for some constant $K$.
Let $r_n=n-\Delta_n$, and $z \in \A(r_n,n)$. There is a constant
$K_b$, independent on $z$ and $n$, such that
\be{seb-main}
\big|\sum_{y\in \A(r_n,n)} G_n(z,y) -\pare{ 
2d\Delta_n\alpha_0(z)-2d(n-\|z\|)^2}\big|\le K_b\pare{(n-\|z\|)\vee 1}\, ,
\ee
with 
\[
\alpha_0(z)=E_z\cro{\|S(H_n)\|-\|z\|\big|H\pare{B^c(0,n)}<H\pare{B(0,r_n)}}.
\]
\ep
\br{rem-d2}
The statement is true in dimension 2, when Green's function
is replaced by the potential kernel
\be{seb-0}
a(x,y)= \E_x \cro{ \sum_{l=0}^{\infty}\id \acc{ S(l)=x}
-\id \acc{ S(l)=y } } \, .
\ee
\er
\bpr
Our strategy is to decompose a path into successive strands lying
entirely in the annulus. The first strand is special since the
starting point is any $z\in \A(r_n,n)$. The other strands, if any,
start all on $\partial \B(0,r_n)$. We estimate the time spent inside
the annulus for each strand. Let us remark that we make use of
three facts: (i) precise asymptotics for Green's function, (ii)
$\{G(0,S(n)),\ n\in \N\}$ is a martingale, and (iii) 
$\{\|S(n)\|^2-n,\ n\in \N\}$ is a martingale.

Choose $z\in \A(r_n,n)$.
We define the following stopping times $(D_i,U_i,\ i\ge 0)$,
corresponding to the ${\mathrm i}^{\mathrm th}$ downward and 
upward crossings of the sphere of radius $r_n$. Let
$\theta(n)$ act on trajectories by time-translation of $n$-units.
Let $\tau=H(B_{r_n})\wedge H_n$, $D_0=U_0=0$, and
\[
D_1=\tau\id_{H(B_{r_n})<H_n}+\infty \id_{H_n<H(B_{r_n})}.
\]
If $D_1<\infty$, then $U_1 =H_{r_n}\circ \theta(D_1)+D_1$, whereas
if $D_1=\infty$, then we set $U_1=\infty$.
We now proceed by induction, and assume $D_i,U_i$ are defined.
If $D_i=\infty$, then $D_{i+1}=\infty$, whereas if $D_i<\infty$, 
(and necessarily $U_i<\infty$) then
\[
D_{i+1}=U_i+\pare{\tau\id_{\tau=H(B_{r_n})}+\infty
\id_{\tau=H_n}}\circ\theta(U_i),\quad\text{and}\quad
U_{i+1} =D_{i+1}+H_{r_n}\circ \theta(D_{i+1}).
\]
With this notation, we can write
\be{tps} 
\begin{split}
\sum_{y\in \A(r_n,n)} G_n(z,y) 
& = E_z\cro{\tau}+
\sum_{i=1}^{\infty} \E_z\cro{ \tau\circ\theta(U_{i})\id_{D_i<\infty}}\\
&=\E_z[\tau]+\P_z\pare{D_1<\infty}\times \I(z),
\end{split}
\ee
where
\be{def-I}
\I(z)=\sum_{i=1}^\infty E_z\cro{ \tau\circ\theta(U_i)\big| D_i<\infty}
\prod_{j=1}^{i-1} \pare{1-P_z(D_{j+1}=\infty\big| D_{j}<\infty)}.
\ee
Now, we compute each term of the right hand side of \reff{tps}.

We have divided the proof in three steps.

\noindent{\bf Step 1:} 
First, we show that there is a positive constant $K$, (independent of
$z$ and $n$) such that when $
z\in \A(r_n,n)$, then
\be{D1bulk}
\big|  \P_z \pare{ D_1<\infty }-
\frac{\alpha_0(z)}{\Delta_n} \big|
 \le \frac{K}{\Delta_n^2}\pare{(n-\|z\|)\vee 1}.
\ee
Note that when $z\in \B(0,n)$, and $n-\|z\|\le 1$, \reff{D1bulk} yields
\be{D1outer}
\big|  \P_z \pare{ D_1<\infty }-
\frac{\E_z\cro{\|S(\tau)\|-\|z\|\big| D_1=\infty}}{\Delta_n} \big|
 \le \frac{K}{\Delta_n^2},
\ee
Secondly, we show that for $z\in \A(r_n,n)$, and $i\ge 1$
\be{D1inner}
\big|  \P_z \pare{ D_{i+1}=\infty\big| D_i<\infty }-
\frac{\E_z\cro{\pare{\|S(U_i)\|-\|S(D_{i+1})}\|
\id_{D_1\circ\theta(U_i)<\infty}
\big| D_{i}<\infty}}{\Delta_n}\big|
\le \frac{K}{\Delta_n^2}.
\ee
Our starting point is the classical Gambler's ruin estimate
\be{gambler-eq}
\P_z \pare{ D_1<\infty } =
\frac{G(0,z)-\E_z \cro{ G(0,S(\tau)) | \, D_1=\infty }}{
\E_z \cro{ G(0,S(\tau)) | \, D_1<\infty }
-\E_z \cro{ G(0,S(\tau)) | \, D_1=\infty }}.
\ee 
We now expand Green's function using asymptotics \reff{main-green}.
For this purpose, it is convenient to define a random variable 
\[
X(z)=\frac{1}{\|z\|}\pare{\|S(\tau)\|^2- \|z\|^2},
\quad\text{ and to set }\eta=\frac{d-2}{2}.
\]
By expressing $S(\tau)$ in terms of $X(z)$, we have
\be{ann-1}
\frac{1}{\|S(\tau)\|^{d-2}}=\frac{1}{\|z\|^{d-2}}
\pare{1+\frac{X(z)}{\|z\|}}^{-\eta}.
\ee
Note that for any $z\in \A(r_n,n)$, $X(z)/\|z\|$ is small. Indeed,
\be{ann-2}
\frac{X(z)}{\|z\|}=\frac{\pare{\|S(\tau)\|-\|z\|}\pare{\|S(\tau)\|+\|z\|}}
{\|z\|^2}
\ee
Since $\Delta_n=n-r_n=O(n^{1/3})$, we have for $n$ large enough
\be{ann-3}
\frac{X(z)}{\|z\|}\le \frac{2(n+1)\Delta_n}{(n-\Delta_n)^2}
\le \frac{8 \Delta_n}{n},\quad
\text{and}\quad 
\sup_{z\in \A(r_n,n)} \pare{ \frac{|X(z)|}{\|z\|}}^3\le 
\frac{8^3 \Delta_n^3}{n}\times \frac{1}{n^2}.
\ee
More precisely, $X(z)$ is of order
$2(\|S(\tau)\|-\|z\|)$. Indeed, 
$\Delta_n^3\le K' n$ for some $K'>0$, and \reff{ann-2} yields
\be{X(z)-general}
X(z)=2\pare{\|S(\tau)\|-\|z\|}+
\pare{\frac{(\|S(\tau)\|-\|z\|)^2}{\|z\|}}\Longrightarrow
\big| X(z)-2\pare{\|S(\tau)\|-\|z\|}\big|\le \frac{K'}{\Delta_n}.
\ee
Finally, we have a constant $K$ such that
\be{ann-4}
\big| \pare{1+\frac{X(z)}{\|z\|}}^{-\eta}-\pare{
1-\eta\frac{X(z)}{\|z\|}+
\eta\frac{\eta+1}{2}\pare{\frac{X(z)}{\|z\|}}^2}\big|
\le \frac{K}{n^2}.
\ee
For any $z\not= 0$, Green's function asymptotics
\reff{main-green} and \reff{ann-4} yields 
\be{ann-5}
\big| G(0,S(\tau))-G(0,z)-\eta C_d\pare{
-\frac{X(z)}{\|z\|^{d-1}}+\frac{\eta+1}{2}\frac{X(z)^2}{\|z\|^{d+1}}}\big|
\le \frac{K}{n^d}.
\ee
Using \reff{gambler-eq} and \reff{ann-5}, we obtain
\be{ann-6}
\P_z \pare{ D_1<\infty }=\frac{ E_z\cro{X(z)|D_1=\infty}-\bar C(z)+
O(\frac{1}{n})}
{ E_z\cro{X(z)|D_1=\infty}- E_z\cro{X(z)|D_1<\infty}+\sous C(z)-
\bar C(z)+O(\frac{1}{n})},
\ee
where
\be{ann-7}
\bar C(z)=\frac{\eta+1}{2} E_z\cro{ \frac{X^2(z)}{\|z\|}\big| D_1=\infty},
\quad\text{and}\quad
\sous C(z)=\frac{\eta+1}{2} E_z\cro{ \frac{X^2(z)}{\|z\|}\big| D_1<\infty}.
\ee
Using \reff{ann-3},
we have some rough estimates on $\bar C$ and $\sous C$. For any
$z\in \A(r_n,n)$, 
\be{C-rough}
\bar C(z)=O\pare{ \frac{\Delta_n^2}{n}}=
O\pare{\frac{1}{\Delta_n}},\quad\text{and}\quad
\sous C(z)=O\pare{ \frac{\Delta_n^2}{n}}=O\pare{\frac{1}{\Delta_n}}.
\ee
Using \reff{X(z)-general}, 
we have better estimates for $\bar C$ and $\sous C$.
\be{C-better}
\bar C(z)=d\frac{(n-\|z\|)^2}{\|z\|}+O\pare{ \frac{\Delta_n}{n}}
,\quad\text{and}\quad
\sous C(z)=d\frac{(\|z\|-r_n)^2}{\|z\|}+O\pare{\frac{\Delta_n}{n}}.
\ee
The rough estimates \reff{C-rough}
allow us to derive from \reff{ann-6}
an estimate for $\P_z(D_1<\infty)$, for any $z\in \A(r_n,n)$.
\be{D1-general}
\begin{split}
\P_z \pare{ D_1<\infty }=&\frac{ E_z\cro{\|S(\tau)\|-\|z\||D_1=\infty}
+O\pare{\frac{1}{\Delta_n}}}
{E_z\cro{\|S(\tau)\|-\|z\||D_1=\infty}-
E_z\cro{\|S(\tau)\|-\|z\||D_1<\infty}+O\pare{\frac{1}{\Delta_n}}}\\
=&\frac{\alpha_0(z)+O\pare{\frac{1}{\Delta_n}}}
{\Delta_n(1+O\pare{\frac{1}{\Delta_n}})}.
\end{split}
\ee
This yields \reff{D1bulk} since $\alpha_0(z)\le
1+(n-\|z\|)\vee 1\le 2 (n-\|z\|)\vee 1$.

\noindent{\bf Case where $z\in \p B(0,r_n)$.}\\
On $\{D_1=\infty\}$, we have
\be{ann-11}
X(z)=2\pare{\|S(\tau)\|-\|z\|}+O(\frac{1}{\Delta_n}).
\ee
On $\{D_1<\infty\}$, we have
\be{ann-12}
X(z)=2\pare{ \|S(\tau)\|-\|z\|}+O(\frac{1}{n}).
\ee
This implies that using \reff{C-better}
\be{C-expansion}
\bar C(z)=d\frac{\Delta_n^2}{\|z\|}+O(\frac{\Delta_n}{n}),
\quad\text{and}\quad \sous C(z)=O(\frac{1}{n}).
\ee
Thus, 
\be{ann-13}
\begin{split}
\P_z\pare{ D_1=\infty}=& \frac{2\E_z\cro{ \|z\|-\|S(\tau)\|
\big| D_1<\infty}+\sous C(z)+O(\frac{1}{n})}
{\E_z\cro{X(z)|D_1=\infty}-\E_z\cro{X(z)|D_1<\infty} +\sous C-\bar C+
O(\frac{1}{n})}\\
=&  \frac{\E_z\cro{ \|z\|-\|S(\tau)\|\big| D_1<\infty}+O(\frac{1}{n})}
{\Delta_n+O(1)}\\
=&\frac{\E_z\cro{ \|z\|-\|S(\tau)\|
\big| D_1<\infty}}{\Delta_n}+O(\frac{1}{\Delta_n^2}).
\end{split}
\ee
In order to obtain \reff{D1inner},
we write \reff{ann-13} on $\{D_i<\infty\}$, and $z=S(U_i)$
as follows. There is a constant $K$ such that on the event
$\{D_{i}<\infty\}$,
\be{seb-bis1}
\big| \E_{S(U_i)}\cro{\id_{ D_{i+1}=\infty}}-
\frac{\E_{S(U_i)}\cro{ \pare{\|S(U_i)\|-\|S(\tau)\|}
\id_{D_1\circ\theta(U_i)<\infty}}}{\Delta_n\times 
P_{S(U_i)}\pare{D_1<\infty}} \big| \le  \frac{K}{\Delta_n^2}.
\ee
Note that \reff{ann-13} implies that $P_{S(U_i)}\pare{D_1<\infty}=1+
O(1/\Delta_n)$, so that \reff{seb-bis1} reads as we integrate over 
$\{D_{i}<\infty\}$ with respect to $E_z$
\be{seb-bis2}
\big| \P_{z}\pare{D_{i+1}=\infty,\ D_{i}<\infty}-
\frac{\E_z\cro{\id_{D_i<\infty}\pare{\|S(U_i)\|-\|S(\tau)\|}
\id_{D_1\circ\theta(U_i)<\infty}}}{\Delta_n}\big| 
\le  \frac{K\P_z(D_i<\infty)}{\Delta_n^2}.
\ee
We obtain \reff{D1inner} as we divide both sides of \reff{seb-bis2}
by $\P_z(D_i<\infty)$.

\noindent{\bf Step 2:}
We show now that for any $z\in \A(r_n,n)$ we have
\be{tau-bulk}
\big| E_z\cro{ \tau}-\pare{d \Delta_n\alpha_0(z)-2d(n-\|z\|)^2}\big|
\le K\pare{(n-\|z\|)\vee 1}.
\ee
When $z\in B_n$ and $n-\|z\|\le 1$,\reff{tau-bulk} reads 
\be{tau-outer}
\big| E_z\cro{ \tau}-\pare{d \Delta_n\alpha_0(z)-2d(n-\|z\|)^2}\le K.
\ee
When $z\in \A(r_n,n)$, and $i\ge 1$, we show that
\be{tau-inner}
\big|\frac{\E_z \cro{ \tau\circ\theta(U_i)\big|D_i<\infty }}{d\Delta_n^2}-
\frac{\E_z\cro{\pare{\|S(U_i)\|-\|S(D_{i+1})\|}\id_{D_1\circ\theta(U_i)
<\infty}\big|D_i<\infty}}{\Delta_n} \big| \le \frac{K}{\Delta_n^2}.
\ee
Using that $\{\|S(n)\|^2-n,\ n\in \N\}$ is a martingale, 
and the optional sampling theorem (see Lemma 2 of \cite{lawler92})
\be{ann-15}
\begin{split}
E_z\cro{\tau}=&\E_z\cro{\|S(\tau)\|^2}-\|z\|^2=\|z\|\times
\E_z\cro{X(z)}\\
=&\|z\|\times \pare{\E_z\cro{X(z)|D_1=\infty}\P_z(D_1=\infty)+
\E_z\cro{X(z)|D_1<\infty}\P_z(D_1<\infty)}.
\end{split}
\ee
Thus, using \reff{ann-6}, simple algebra yields 
\be{ann-16}
\E_z\cro{\tau}=\|z\|\times \pare{(\sous C(z)-\bar C(z))
\P_z(D_1<\infty)+\bar C(z)}+O(1).
\ee
By recalling \reff{C-better} and \reff{D1bulk}
\be{ann-17}
\begin{split}
\E_z\cro{\tau}=&d \pare{ \pare{ (\|z\|-r_n)^2-(n-\|z\|)^2+O(\Delta_n)}
\pare{ \frac{\alpha_0(z)}{\Delta_n}+O(\frac{(n-\|z\|)\vee 1}{\Delta_n^2})}}+O(1)\\
=&d(2\|z\|-n-r_n)\alpha_0(z)+O\pare{(n-\|z\|)\vee 1}\\
=&d\Delta_n \alpha_0(z)-2d(n-\|z\|)^2+O\pare{(n-\|z\|)\vee 1}
\end{split}
\ee
Note that in the case where $n-\|z\|\le 1$, \reff{ann-17} yields
\reff{tau-outer}.

Assume now that $z\in \p B(0,r_n)$. From \reff{ann-16}, we have
\be{ann-19}
\E_z\cro{\tau}=\|z\|\times \pare{(\bar C(z)-\sous C(z))
\P_z(D_1=\infty)+\sous C(z)}+O(1).
\ee
We use \reff{D1inner}, \reff{ann-11} and \reff{ann-12} to obtain
\be{ann-20}
\begin{split}
\E_z\cro{\tau}=& \|z\|\pare{\pare{ d\Delta_n^2+O(\Delta_n)}
\pare{\frac{\E_z\cro{ \|z\|-\|S(\tau)\|\big| D_1<\infty}}{\Delta_n}
+O(\frac{1}{\Delta_n^2})}}+O(1)\\
=& d\Delta_n \E_z\cro{ \|z\|-\|S(\tau)\|\big| D_1<\infty}+O(1).
\end{split}
\ee
Now, write \reff{ann-20} as follows. There is a constant $K$ such that
for any $z\in \p \B(0,r_n)$
\be{ann-21}
\big| \frac{\E_z\cro{\tau}}{d\Delta_n^2}-
\frac{\E_z\cro{\|z\|-\|S(\tau)\|\id_{D_1<\infty}}}{\Delta_n \P_z(D_1<\infty)}
\big|\le \frac{K}{\Delta_n^2}.
\ee
Note that by \reff{ann-13}, we have that
$\Delta_n \P_z(D_1<\infty)=\Delta_n+O(1)$, and $|\|z\|-\|S(\tau)\|
\id_{D_1<\infty}|\le 1$, we have
\be{ann-main}
\big| \frac{\E_z\cro{\tau}}{d\Delta_n^2}-
\frac{\E_z\cro{\|z\|-\|S(\tau)\|\id_{D_1<\infty}}}{\Delta_n}
\big|\le \frac{K}{\Delta_n^2}.
\ee
We replace $z$ by $S(U_i)$ in \reff{ann-main} under the event $\{D_i<\infty
\}$ to obtain
\be{ann-22}
\big| \frac{\E_{S(U_i)}\cro{\tau}}{d\Delta_n^2}-
\frac{\E_{S(U_i)}\cro{
\pare{\|S(U_i)\|-\|S(D_1\circ\theta(U_i))\|}
\id_{D_1\circ\theta(U_i)<\infty}}}
{\Delta_n} \big|\le \frac{K}{\Delta_n^2},
\ee
We multiply both sides of \reff{ann-22} by $\id_{D_i<\infty}$,
take the expectation on both side of \reff{ann-22}, and divide
by $\P_z(D_i<\infty)$ to obtain \reff{tau-inner}.

\noindent{\bf Step 3:} For $i\ge 1$, we show the following bounds 
\be{ai}
2\ge \gamma_i\geq \frac{1}{4d\sqrt{d}},
\quad\text{where }
\gamma_i=\E_{z}
\cro{\pare{\|S(U_i)\|-\|S(D_{i+1})\|}\id_{D_{i+1}<\infty}
\big| D_i<\infty}.
\ee
The upper bound is obvious. For the lower bound, first we restrict to
$\{D_i<\infty\}$, so that $U_i<\infty$.
By Lemma~\ref{lem-border}, $S(U_i)$ has a nearest neighbor $x$,
within $\B(0,r_n)$ such that 
$\norm{S(U_i)}- \norm{x} \geq 1/(2\sqrt{d})$,
and \reff{ai} is immediate.

\noindent{\bf Step 4:} We show \reff{seb-main} using \reff{tps}.
For $p$ such that $1\leq p \leq \infty$, let
\be{seb-9}
\sigma_p=\sum_{i=1}^{p} 
\E_z \cro{\tau\circ\theta(U_i)\big| D_i<\infty}\\
\prod_{j=1}^{i-1} \pare{1-\P_z\pare{D_{j+1}=\infty\big|D_j<\infty}}.
\ee
Now, \reff{def-I} reads
$\I(z)=\lim_{p\to\infty} \sigma_p$.
We establish in this step that, for some constant $\tilde K$, any
integer $n$
\be{step4}
\lim_{p\to\infty}|1-\frac{\sigma_p}{d\Delta_n^2}|\le \frac{\tilde K}{\Delta_n}.
\ee
Once we prove \reff{step4}, we have all the bounds to
estimate the right hand side of \reff{tps}. Indeed, using
\reff{tau-bulk}, \reff{D1bulk} and \reff{step4}, we have
\be{conclu-appendix}
\begin{split}
\E_z[\tau]+\P_z\pare{D_1<\infty}\times \I(z)=&
d\Delta_n\alpha_0(z)-2d(n-\|z\|)^2 +O\pare{(n-\|z\|)\vee 1} \\
&\ +\pare{\frac{\alpha_0(z)}{\Delta_n}
+O\pare{\frac{(n-\|z\|)\vee 1}{\Delta_n^2}}}
\times\pare{d \Delta_n^2+O(\Delta_n)}\\
&=2d\Delta_n\alpha_0(z)-2d\pare{n-\|z\|}^2+O\pare{(n-\|z\|)\vee 1}.
\end{split}
\ee
In order now to prove \reff{step4}, we
introduce first some shorthand notation. For $p$ and $j$ positive integers
\be{def-amin}
a_p=1-\frac{\sigma_p}{d\Delta_n^2},\quad
\alpha_j=P_z(D_{j+1}=\infty|D_{j}<\infty),\quad\text{and}\quad
\beta_j= \frac{\E_z \cro{\tau\circ\theta(U_j)\big| D_j<\infty}}
{d\Delta_n^2}.
\ee
With this notation \reff{D1inner} and \reff{tau-inner} read as follows.
\be{step2}
|\alpha_j-\frac{\gamma_j}{\Delta_n}|\le \frac{K}{\Delta^2_n},\quad
\text{and}\quad
|\beta_j-\frac{\gamma_j}{\Delta_n}|\le \frac{K}{\Delta^2_n},\quad
\text{so that}\quad
|\alpha_j-\beta_j|\le \frac{2K}{\Delta^2_n}.
\ee
Let us rewrite \reff{seb-9} as
\be{amin-recursion}
a_p=a_{p-1}-\beta_p \prod_{j=1}^{p-1} (1-\alpha_j).
\ee
In order to establish \reff{step4}, we show by induction that
\be{seb-induction}
\big| a_p-\prod_{j=1}^{p} (1-\alpha_j)|\le \epsilon_p,
\ee
with for $p>1$
\be{seb16}
\epsilon_p=
\epsilon_{p-1}+\frac{2K}{\Delta_n^2} \prod_{j=1}^{p-1} (1-\alpha_j)
\quad\text{and}\quad \epsilon_1=\frac{2K}{\Delta_n^2}.
\ee
Note that it is easy to estimate $\epsilon_p$ from \reff{seb16}.
There is a constant $\kappa_S$ such that
\[
\epsilon_d\le\frac{2K}{\Delta_n^2}(1+\sum_{k=1}^p\exp(-\sum_{j=1}^k 
\alpha_j))
\le  \frac{2K}{\Delta_n^2}\pare{1+\sum_{k=1}^p\exp\pare{-\sum_{j=1}^k 
\frac{\gamma_j}{2\Delta_n}}} \le\frac{2K}{\Delta^2_n}\kappa_S\Delta_n=
\frac{2K\kappa_S}{\Delta_n}.
\]
Now, \reff{seb-induction} holds for $p=1$, and assume it holds for $p-1$.
Then
\be{seb-10}
(1-\beta_p) \prod_{j=1}^{p-1} (1-\alpha_j)-\epsilon_{p-1}\le a_p\le
(1-\beta_p) \prod_{j=1}^{p-1} (1-\alpha_j)+\epsilon_{p-1}.
\ee
Then by \reff{step2}, we have \reff{seb-induction} with $\epsilon_p$
satisfying \reff{seb16}.

Now \reff{step4} follows as we notice that Step 3 implies that
\[
\lim_{p\to\infty} \prod_{j=1}^{p} (1-\alpha_j)=0.
\]
\epr

\noindent{\bf Acknowledgements.}
The authors thank the CIRM for a friendly atmosphere
during their stay as part of the research in pairs program.

\end{document}